# X-PIPS-M DATA SUMMARY

Progress through Calculus: 2017-19

## Abstract

A compilation of descriptive statistics of responses to the X-PIPS-M suite of surveys for students, student instructors, and faculty instructors of introductory mathematics courses. Data were obtained from Precalculus, Calculus 1, and Calculus 2 courses at 12 universities during the 2017-18 and 2018-19 academic years, as part of the Progress through Calculus research project (NSF DUE IUSE #1430540).

Ciera Street, Naneh Apkarian, Jessica Gehrtz, Rachel Tremaine, Victoria Barron, Matthew Voigt, & Jessica Ellis Hagman on behalf of the Progress through Calculus Research Team.

## Table of Contents

*This material is based upon work supported by the National Science Foundation under grant DUE IUSE #1430540. Any opinions, findings, and conclusions or recommendations expressed in this material are those of the authors and project team and do not necessarily reflect the views of the National Science Foundation.*

# Overview / Introduction

## Project Overview

Data for this report come from the *Progress through Calculus* (PtC) project ([maa.org/PtC](maa.org/PtC)), which was sponsored by the National Science Foundation (DUE IUSE #1430540) and run with support from the Mathematical Association of America. The project focused on the Precalculus to Calculus 2 (P2C2) course sequence, which is the course sequence generally required of undergraduate STEM majors. This multiyear project, begun in 2015, consisted of two phases. Phase 1 involved a survey of all university mathematics departments in the United States. The aim of this phase was to catalogue common practices, and in particular, to document the prevalence of program characteristics previously shown to support student success in Calculus 1 ([Bressoud, Mesa, & Rasmussen, 2015](#)). Results of this first phase have been published as a journal article ([Rasmussen et al., 2019](#)) and in a technical report like this one ([Apkarian & Kirin, 2017](#)). Phase 2 of the PtC project consisted of multi-year case studies of 12 selected mathematics programs. This report summarizes responses to surveys collected as part of Phase 2.

## Survey Design & Development

The surveys reported on here were co-developed with members of the SEMINAL project ([aplu.org/SEMINAL](aplu.org/SEMINAL)), and the text of these questions, along with a report of the development process, is available as a white paper ([Apkarian et al., 2019](#)). Given this availability, we omit a lengthy discussion of the survey development in this technical summary report. In brief, the X-PIPS-M refers to a suite of instruments created from postsecondary instructional practice surveys (PIPS) about mathematics (the *M*) for different stakeholders (the *X*). The core of these surveys draws on the *Postsecondary Instructional Practices Survey* for STEM instructors which was developed elsewhere (Walter et al., 2016). Our version amended the original items slightly for the mathematical context and contributed additional questions about context and attitudes. We also developed versions specifically for students, faculty instructors, and student instructors to capture different aspects of introductory mathematics education. Surveys generally consisted of (1) questions about the participant, including demographic and contextual questions; and (2) questions about course delivery and experience.

## Data Collection & Participants

Data were collected at 12 universities in each term of the 2017-18 academic year and the Fall term of 2019. Surveys were distributed to all instructors of, student instructors of, and students enrolled in mainstream P2C2 courses. Graduate students and undergraduate students who served as either instructor of record or teaching assistant were considered to be student instructors. For our purposes, "Precalculus" refers to courses which function as immediate prerequisites for single-variable calculus regardless of their official title; "Calculus 1" refers to a first course in single-variable calculus, and "Calculus 2" refers to subsequent courses in single variable calculus. The term "mainstream" indicates that a course functions as part of the prerequisite sequence for additional upper-division mathematics courses; terminal courses that incorporate similar content were not included. Surveys were appropriately tailored to each context: students saw the exact name of their course throughout (indicated in this document as [course]) and students enrolled in courses without an accompanying lab/recitation component were not presented with questions about labs/recitations. Select items were presented based on responses to a previous question; this is indicated where relevant.

Surveys were distributed via email to all participants and the survey was completed electronically using Qualtrics. The research team contacted instructors and student instructors directly, while students were contacted by their instructors; this distribution was facilitated by course coordinators and program liaisons as fit the local context. Instructors were encouraged to offer extra credit to their students for accessing the survey, but this was not required and not all instructors offered such incentives.

Due to the study design and the nature of P2C2 programs, many participants appear more than once in the data. People are double counted if they responded multiple times, which might occur for (1) students across multiple terms; (2) instructors teaching more than one P2C2 course in a single term; (3) instructors teaching any P2C2 course across multiple terms in the study. For instructors who reported on more than one course they were teaching in a particular term, their demographic/contextual information is reported only once while each course-related response is reported separately. For people who responded to the survey multiple times for the same course in the same term, we marked these as duplicates and only report the responses from the most complete survey attempt for each participant.

## Reading this Report

This report proceeds in three parts. The student survey (S-PIPS) is presented first, starting on page 5; student instructor survey (U-GPIPS) results begin on page 29; faculty instructors (PIPS) results begin on page 43. Within each, we report first on participant characteristics (demographics, identity, and context) and then on responses to the questions about course delivery and experiences. Participant characteristics are aggregated, while course specific responses are separated into responses pertaining to Precalculus, Calculus 1, and Calculus 2.

Some response options and Likert statements have been abbreviated in tables for readability. This is made clear in the text describing the survey item, and when questions or statements have been abbreviated, they are indicated with an asterisk (*). Exact wording of all items is available in the aforementioned white paper (Apkarian et al., 2019).

Depending on the survey question, tabulated data are reported in one of two ways to show both the numeric count and the proportion of respondents who selected particular choices. In condensed tables, responses are shown as *number (proportion)*; in extended tables a row of counts appears directly above a row of proportions. Proportions are either out of the total who responded to a particular question (single-response items) or out of the grand total (multiple-response items); any deviations from this are noted in the text.

# S-PIPS-M: Student Survey Results

The S-PIPS data contain responses from 19,191 undergraduate students from 12 universities across the nation. We first present the self-identified demographic and contextual information of these students, followed by students' reports of their experiences in their P2C2 courses.

## Student Identity & Individual Context

These questions asked about students' identity, including demographic information and academic context. Questions about Gender, Race and Ethnicity, and Sexual Orientation were presented as "select all that apply" items with the option to self-identify descriptors not on our lists. We intentionally provided more nuanced options than are generally used in demographic reports (e.g., the US census) as we value the complex nature of identity. Responses to these items were recoded for interpretability in this report; details are provided alongside each item.

### Gender

Do you consider yourself to be (Select all that apply):

| Cisgender Man | Cisgender Woman | Gender Non-Binary | Transgender | Transgender Man | Transgender Woman | Not reported |
|---|---|---|---|---|---|---|
| 10045 | 7972 | 39 | 48 | 14 | 9 | 1064 |
| 0.523 | 0.415 | 0.002 | 0.003 | 0.001 | 0.000 | 0.055 |

*Survey response options: Man, Transgender, Woman, Not listed (please specify), Prefer not to disclose.*

*Selections of only Man, only Woman, or only Transgender are coded as Cisgender Man, Cisgender Woman, and Transgender respectively; selections of Transgender and Man or Transgender and Woman are coded as Transgender Man and Transgender Woman, respectively; selections of Man and Woman (with or without additional selections) or self-described non-binary gender identities (e.g., genderfluid, agender) are coded as Gender Non-Binary; selections of Prefer not to disclose and no selection were coded as Not reported.*

### Race and Ethnicity

Do you consider yourself to be (Select all that apply):

| | |
|---|---|
| Asian | 1449 (0.076) |
| Alaska Native or Native American | 70 (0.004) |
| Black or African American | 1474 (0.077) |
| Hispanic or Latinx | 2137 (0.111) |
| Middle Eastern or North African | 648 (0.034) |
| Multiple Race/Ethnicity | 2044 (0.107) |
| Native Hawaiian or Pacific Islander | 42 (0.002) |
| White | 9962 (0.519) |
| Not reported | 1365 (0.071) |

*Survey response options: Alaskan Native or Native American, Black or African American, Central Asian, East Asian, Hispanic or Latinx, Middle Eastern or North African, Native Hawaiian or Pacific Islander, Southeast Asian, South Asian, White, Not listed (please specify), Prefer not to disclose.*

*Selections of more than one of Central Asian, East Asian, Southeast Asian, South Asian without any other choices are coded as Asian; selections of multiple responses (aside from the previous) are coded as Multiple Race/Ethnicity; selections of Prefer not to disclose and no selection are coded as Not reported.*

## Sexual Orientation

Do you consider yourself to be (Select all that apply):

| Straight | Gay | Lesbian | Bisexual | Asexual | Queer | Queer+ | Straight+ | Not reported |
|---|---|---|---|---|---|---|---|---|
| 15618 | 234 | 109 | 698 | 495 | 73 | 169 | 115 | 1680 |
| 0.814 | 0.012 | 0.006 | 0.036 | 0.026 | 0.004 | 0.009 | 0.006 | 0.088 |

*Survey response options: Asexual, Bisexual, Gay, Lesbian, Queer, Straight (Heterosexual), Not listed (please specify), Prefer not to disclose.*

*Categories are based on research from Voigt (2020) where selections of exactly one listed option are not recoded; selections of Straight (Heterosexual) and one or more additional options are coded as Straight+; selections of more than one option but not Straight (Heterosexual) are coded as Queer+; selections of Prefer not to disclose and no selection are coded as Not reported.*

## Additional Context / Identity Markers

Do you consider yourself to be (Select all that apply):

| | |
|---|---|
| International Student | 1507 (0.079) |
| First-generation college student | 4963 (0.259) |
| Commuter student | 3234 (0.169) |
| Transfer student | 1927 (0.1) |
| Student with a disability | 558 (0.029) |
| Student athlete | 973 (0.051) |
| Current or former English language learner | 1003 (0.052) |
| Parent or guardian | 135 (0.007) |

## FAFSA and Pell Grants

Did you use FAFSA to apply for financial aid?

| Yes | No | Not reported |
|---|---|---|
| 13,164 (0.686) | 4102 (0.214) | 1925 (0.1) |

Did you receive a free grant (e.g., Pell Grant)? [Shown only to the 13,164 who responded "Yes" to the previous].

| Yes | No | I don't know | Not reported |
|---|---|---|---|
| 5600 (0.425) | 4933 (0.375) | 2086 (0.159) | 545 (0.041) |

## Work Hours Per Week

Approximately how many hours per week did you work at a job this term?

| 0 | 1 – 5 | 6 – 10 | 11 – 15 | 16 – 20 | 21 – 30 | >30 | Not reported |
|---|---|---|---|---|---|---|---|
| 9943 | 1067 | 1417 | 1224 | 1516 | 1309 | 1013 | 884 |
| 0.541 | 0.058 | 0.077 | 0.067 | 0.083 | 0.071 | 0.055 | 0.048 |

## Class Standing

What is your class standing?

| First Year | Sophomore | Junior | Senior | Other | Not reported |
|---|---|---|---|---|---|
| 10359 | 4885 | 2113 | 678 | 209 | 176 |
| 0.562 | 0.265 | 0.115 | 0.037 | 0.011 | 0.01 |

## STEM Intention

Have you declared, or do you intend to declare, a STEM (science, technology, engineering, or mathematics) major?

| Yes | No | Unsure | Not reported |
|---|---|---|---|
| 13121 | 3386 | 1520 | 364 |
| 0.713 | 0.184 | 0.083 | 0.02 |

## Math Preparation

Do you think your previous math courses adequately prepared you for [current course]?

| Yes | No |
|---|---|
| 15010 | 3352 |
| 0.817 | 0.183 |

## Student Experiences

These questions of the survey relate to the students' reports of their experiences in their Precalculus, Calculus I, and Calculus II courses. The respondents totaled 7136 in Precalculus, 6888 in Calculus I, and 7136 in Calculus II.

### Percentages of Class Time

What percent of regular class time, over the whole term, did you spend… [must total 100]

|  | n | mean | sd | median |
| --- | --- | --- | --- | --- |
| **All Students** | | | | |
| *Working on tasks individually* | 19063 | 23.22 | 23.08 | 17 |
| *Working on tasks in small groups* | 19063 | 11.77 | 16.25 | 5 |
| *Engaging in whole class discussion* | 19063 | 14.49 | 16.26 | 10 |
| *Listening to the instructor lecture or solve problems* | 19063 | 50.52 | 28.25 | 50 |
| **Precalculus Students** | | | | |
| *Working on tasks individually* | 7085 | 28.51 | 25.55 | 24 |
| *Working on tasks in small groups* | 7085 | 13.13 | 18.46 | 5 |
| *Engaging in whole class discussion* | 7085 | 15.46 | 16.88 | 10 |
| *Listening to the instructor lecture or solve problems* | 7085 | 42.9 | 27.33 | 41 |
| **Calculus 1 Students** | | | | |
| *Working on tasks individually* | 6819 | 21.41 | 21.96 | 15 |
| *Working on tasks in small groups* | 6819 | 9.963 | 14.55 | 3 |
| *Engaging in whole class discussion* | 6819 | 13.56 | 15.74 | 10 |
| *Listening to the instructor lecture or solve problems* | 6819 | 55.07 | 28.13 | 52 |
| **Calculus 2 Students** | | | | |
| *Working on tasks individually* | 5159 | 18.36 | 19.19 | 13 |
| *Working on tasks in small groups* | 5159 | 12.27 | 14.87 | 8 |
| *Engaging in whole class discussion* | 5159 | 14.4 | 16 | 10 |
| *Listening to the instructor lecture or solve problems* | 5159 | 54.97 | 27.41 | 52 |

## Percentages of Recitation/Lab Time

What percent of regular **recitation/lab time**, over the whole term, did you spend… [must total 100]

|  | n | mean | sd | median |
|---|---|---|---|---|
| **All Students** | | | | |
| *Working on tasks individually* | 7040 | 29.94 | 28.84 | 24 |
| *Working on tasks in small groups* | 7040 | 41.19 | 32.82 | 40 |
| *Engaging in whole class discussion* | 7040 | 9.82 | 14.11 | 2 |
| *Listening to the instructor lecture or solve problems* | 7040 | 19.05 | 23.07 | 10 |
| **Precalculus Students** | | | | |
| *Working on tasks individually* | 1583 | 33.55 | 32.61 | 25 |
| *Working on tasks in small groups* | 1583 | 48.95 | 32.21 | 50 |
| *Engaging in whole class discussion* | 1583 | 7.196 | 12.43 | 0 |
| *Listening to the instructor lecture or solve problems* | 1583 | 10.3 | 16.21 | 0 |
| **Calculus 1 Students** | | | | |
| *Working on tasks individually* | 3929 | 30.11 | 27.7 | 25 |
| *Working on tasks in small groups* | 3929 | 38.57 | 32.22 | 35 |
| *Engaging in whole class discussion* | 3929 | 10.51 | 14.3 | 5 |
| *Listening to the instructor lecture or solve problems* | 3929 | 20.81 | 23.26 | 13 |
| **Calculus 2 Students** | | | | |
| *Working on tasks individually* | 1528 | 25.8 | 26.98 | 20 |
| *Working on tasks in small groups* | 1528 | 39.87 | 33.76 | 35 |
| *Engaging in whole class discussion* | 1528 | 10.76 | 14.91 | 4.5 |
| *Listening to the instructor lecture or solve problems* | 1528 | 23.57 | 26.09 | 15 |

## Overall Experience

Indicate the degree to which the following statements describe your experience in [course]:

| | | Never | Minimally | Somewhat | Mostly | Very |
|---|---|---|---|---|---|---|
| The test questions focus on important facts and definitions from the course. | PC | 267 | 562 | 1517 | 2346 | 2399 |
| | | 0.038 | 0.079 | 0.214 | 0.331 | 0.338 |
| | C1 | 135 | 541 | 1438 | 2623 | 2121 |
| | | 0.02 | 0.079 | 0.21 | 0.382 | 0.309 |
| | C2 | 102 | 337 | 973 | 1854 | 1876 |
| | | 0.02 | 0.066 | 0.189 | 0.361 | 0.365 |

|  |  | Never | Minimally | Somewhat | Mostly | Very |
|---|---|---|---|---|---|---|
| The test questions require me to apply course concepts to unfamiliar situations. | PC | 783 | 1256 | 2013 | 1814 | 1231 |
|  |  | 0.11 | 0.177 | 0.284 | 0.256 | 0.173 |
|  | C1 | 443 | 1092 | 1896 | 1997 | 1433 |
|  |  | 0.065 | 0.159 | 0.276 | 0.291 | 0.209 |
|  | C2 | 402 | 897 | 1466 | 1430 | 947 |
|  |  | 0.078 | 0.174 | 0.285 | 0.278 | 0.184 |
| I use technology or online resources in relation to this course. | PC | 259 | 422 | 1086 | 1862 | 3476 |
|  |  | 0.036 | 0.059 | 0.153 | 0.262 | 0.489 |
|  | C1 | 276 | 534 | 1171 | 2015 | 2872 |
|  |  | 0.04 | 0.078 | 0.171 | 0.293 | 0.418 |
|  | C2 | 226 | 418 | 903 | 1509 | 2088 |
|  |  | 0.044 | 0.081 | 0.176 | 0.293 | 0.406 |
| I make connections between related ideas or concepts when completing assignments. | PC | 402 | 714 | 2020 | 2281 | 1688 |
|  |  | 0.057 | 0.1 | 0.284 | 0.321 | 0.238 |
|  | C1 | 236 | 733 | 1977 | 2497 | 1424 |
|  |  | 0.034 | 0.107 | 0.288 | 0.364 | 0.207 |
|  | C2 | 128 | 465 | 1480 | 1925 | 1148 |
|  |  | 0.025 | 0.09 | 0.288 | 0.374 | 0.223 |
| I receive feedback on my assignments without being assigned a formal grade. | PC | 2540 | 1473 | 1500 | 873 | 712 |
|  |  | 0.358 | 0.208 | 0.211 | 0.123 | 0.1 |
|  | C1 | 2675 | 1635 | 1372 | 736 | 444 |
|  |  | 0.39 | 0.238 | 0.2 | 0.107 | 0.065 |
|  | C2 | 1948 | 1292 | 1044 | 542 | 319 |
|  |  | 0.379 | 0.251 | 0.203 | 0.105 | 0.062 |
| I see my instructor(s) outside of class for help. | PC | 3208 | 1563 | 1150 | 593 | 591 |
|  |  | 0.452 | 0.22 | 0.162 | 0.083 | 0.083 |
|  | C1 | 2757 | 1751 | 1106 | 612 | 630 |
|  |  | 0.402 | 0.255 | 0.161 | 0.089 | 0.092 |
|  | C2 | 1759 | 1330 | 929 | 591 | 531 |
|  |  | 0.342 | 0.259 | 0.181 | 0.115 | 0.103 |

|  |  | Never | Minimally | Somewhat | Mostly | Very |
|---|---|---|---|---|---|---|
| I work with peers outside of class on math problems. | PC | 2168 | 1211 | 1367 | 1186 | 1172 |
|  |  | 0.305 | 0.17 | 0.192 | 0.167 | 0.165 |
|  | C1 | 1527 | 1258 | 1456 | 1288 | 1340 |
|  |  | 0.222 | 0.183 | 0.212 | 0.188 | 0.195 |
|  | C2 | 1055 | 962 | 1066 | 1011 | 1049 |
|  |  | 0.205 | 0.187 | 0.207 | 0.197 | 0.204 |
| I attend tutoring sessions outside of class time. | PC | 3030 | 1260 | 1091 | 754 | 972 |
|  |  | 0.426 | 0.177 | 0.154 | 0.106 | 0.137 |
|  | C1 | 2838 | 1291 | 1071 | 759 | 904 |
|  |  | 0.414 | 0.188 | 0.156 | 0.111 | 0.132 |
|  | C2 | 2043 | 1014 | 805 | 623 | 661 |
|  |  | 0.397 | 0.197 | 0.156 | 0.121 | 0.128 |

## Technology Use

Which technologies and/or online resources do you use? Mark all that apply.

|  |  |  |
|---|---|---|
| Graphing Calculator | PC | 3511 (0.492) |
|  | C1 | 3491 (0.507) |
|  | C2 | 2376 (0.46) |
| Clickers or other polling devices | PC | 324 (0.045) |
|  | C1 | 210 (0.03) |
|  | C2 | 188 (0.036) |
| Computer algebra software (e.g. Maple, Matlab, Mathematica) | PC | 1965 (0.275) |
|  | C1 | 1371 (0.199) |
|  | C2 | 1129 (0.219) |
| Online search engines (e.g. Google) | PC | 3698 (0.518) |
|  | C1 | 4066 (0.59) |
|  | C2 | 3090 (0.598) |
| Online textbooks | PC | 3059 (0.429) |
|  | C1 | 2893 (0.42) |
|  | C2 | 1992 (0.386) |

| Online tutorials (e.g. Khan Academy, YouTube videos) | PC | 4575 (0.641) |
|---|---|---|
| | C1 | 4740 (0.688) |
| | C2 | 3633 (0.703) |
| Online computational or graphing tools (e.g. Wolfram Alpha, Geogebra, Desmos) | PC | 2593 (0.363) |
| | C1 | 3420 (0.497) |
| | C2 | 2972 (0.575) |
| Online homework (e.g. WebAssign, MyMathLab, WebWork) | PC | 5072 (0.711) |
| | C1 | 4344 (0.631) |
| | C2 | 2551 (0.494) |
| Online forums (e.g. Chegg, StackExchange, Slader) | PC | 892 (0.125) |
| | C1 | 1623 (0.236) |
| | C2 | 1480 (0.286) |
| Learning management systems (e.g. Blackboard, Canvas, Piazza) | PC | 1873 (0.262) |
| | C1 | 1798 (0.261) |
| | C2 | 1419 (0.275) |

## Classroom Experience

Indicate the degree to which the following statements describe your experience in regular course meetings of [course]:

| | | Not at all | Minimally | Somewhat | Mostly | Very |
|---|---|---|---|---|---|---|
| I listen as the instructor guides me through major topics. | PC | 169 | 262 | 919 | 2123 | 3558 |
| | | 0.024 | 0.037 | 0.131 | 0.302 | 0.506 |
| | C1 | 114 | 202 | 871 | 2129 | 3486 |
| | | 0.017 | 0.03 | 0.128 | 0.313 | 0.512 |
| | C2 | 36 | 119 | 434 | 1578 | 2925 |
| | | 0.007 | 0.023 | 0.085 | 0.31 | 0.574 |
| The class activities connect course content to my life and future work. | PC | 1495 | 1458 | 1875 | 1265 | 939 |
| | | 0.213 | 0.207 | 0.267 | 0.18 | 0.134 |
| | C1 | 1107 | 1535 | 1961 | 1331 | 861 |
| | | 0.163 | 0.226 | 0.289 | 0.196 | 0.127 |
| | C2 | 734 | 1078 | 1465 | 1118 | 699 |
| | | 0.144 | 0.212 | 0.288 | 0.219 | 0.137 |

|  |  | Not at all | Minimally | Somewhat | Mostly | Very |
|---|---|---|---|---|---|---|
| I receive immediate feedback on my work during class (e.g., clickers). | PC | 2017 | 1183 | 1392 | 1203 | 1243 |
|  |  | 0.287 | 0.168 | 0.198 | 0.171 | 0.177 |
|  | C1 | 2504 | 1291 | 1296 | 956 | 757 |
|  |  | 0.368 | 0.19 | 0.19 | 0.141 | 0.111 |
|  | C2 | 1662 | 975 | 953 | 853 | 645 |
|  |  | 0.327 | 0.192 | 0.187 | 0.168 | 0.127 |
| I am asked to respond to questions during class time. | PC | 1357 | 1280 | 1710 | 1432 | 1257 |
|  |  | 0.193 | 0.182 | 0.243 | 0.204 | 0.179 |
|  | C1 | 1525 | 1485 | 1676 | 1271 | 839 |
|  |  | 0.224 | 0.219 | 0.247 | 0.187 | 0.123 |
|  | C2 | 857 | 1089 | 1294 | 1077 | 775 |
|  |  | 0.168 | 0.214 | 0.254 | 0.212 | 0.152 |
| *In my class a variety of means (graphs, symbols, tables, etc.) are used to represent course topics and/or solve problems. | PC | 356 | 683 | 1645 | 2038 | 2314 |
|  |  | 0.051 | 0.097 | 0.234 | 0.29 | 0.329 |
|  | C1 | 269 | 707 | 1655 | 2161 | 2008 |
|  |  | 0.04 | 0.104 | 0.243 | 0.318 | 0.295 |
|  | C2 | 210 | 547 | 1223 | 1615 | 1505 |
|  |  | 0.041 | 0.107 | 0.24 | 0.317 | 0.295 |
| I talk with other students about course topics during class. | PC | 1302 | 1190 | 1586 | 1566 | 1383 |
|  |  | 0.185 | 0.169 | 0.226 | 0.223 | 0.197 |
|  | C1 | 1297 | 1309 | 1665 | 1496 | 1018 |
|  |  | 0.191 | 0.193 | 0.245 | 0.22 | 0.15 |
|  | C2 | 858 | 961 | 1174 | 1171 | 931 |
|  |  | 0.168 | 0.189 | 0.23 | 0.23 | 0.183 |
| I constructively criticize other students' ideas during class. | PC | 3388 | 1424 | 1172 | 615 | 426 |
|  |  | 0.482 | 0.203 | 0.167 | 0.088 | 0.061 |
|  | C1 | 3351 | 1579 | 1009 | 545 | 318 |
|  |  | 0.493 | 0.232 | 0.148 | 0.08 | 0.047 |
|  | C2 | 2233 | 1258 | 846 | 501 | 258 |
|  |  | 0.438 | 0.247 | 0.166 | 0.098 | 0.051 |

|  |  | Not at all | Minimally | Somewhat | Mostly | Very |
|---|---|---|---|---|---|---|
| I discuss the difficulties I have with math with other students during class. | PC | 1340 | 1114 | 1573 | 1501 | 1506 |
|  |  | 0.191 | 0.158 | 0.224 | 0.213 | 0.214 |
|  | C1 | 1265 | 1281 | 1574 | 1492 | 1188 |
|  |  | 0.186 | 0.188 | 0.231 | 0.219 | 0.175 |
|  | C2 | 902 | 901 | 1182 | 1186 | 921 |
|  |  | 0.177 | 0.177 | 0.232 | 0.233 | 0.181 |
| I work on problems individually during class time. | PC | 629 | 1090 | 1882 | 1735 | 1700 |
|  |  | 0.089 | 0.155 | 0.267 | 0.247 | 0.242 |
|  | C1 | 821 | 1486 | 1830 | 1496 | 1165 |
|  |  | 0.121 | 0.219 | 0.269 | 0.22 | 0.171 |
|  | C2 | 634 | 1169 | 1434 | 1098 | 754 |
|  |  | 0.125 | 0.23 | 0.282 | 0.216 | 0.148 |
| I work with other students in a small group during class. | PC | 2112 | 1248 | 1272 | 1084 | 1320 |
|  |  | 0.3 | 0.177 | 0.181 | 0.154 | 0.188 |
|  | C1 | 2314 | 1427 | 1163 | 942 | 948 |
|  |  | 0.341 | 0.21 | 0.171 | 0.139 | 0.14 |
|  | C2 | 1439 | 942 | 975 | 834 | 899 |
|  |  | 0.283 | 0.185 | 0.192 | 0.164 | 0.177 |
| Multiple approaches to solving a problem are discussed in class. | PC | 429 | 820 | 1713 | 2061 | 2001 |
|  |  | 0.061 | 0.117 | 0.244 | 0.293 | 0.285 |
|  | C1 | 359 | 902 | 1845 | 2113 | 1565 |
|  |  | 0.053 | 0.133 | 0.272 | 0.311 | 0.231 |
|  | C2 | 170 | 501 | 1247 | 1681 | 1494 |
|  |  | 0.033 | 0.098 | 0.245 | 0.33 | 0.293 |
| I have enough time during class to reflect about the processes I use to solve problems. | PC | 715 | 1176 | 1873 | 1815 | 1455 |
|  |  | 0.102 | 0.167 | 0.266 | 0.258 | 0.207 |
|  | C1 | 706 | 1344 | 1852 | 1797 | 1099 |
|  |  | 0.104 | 0.198 | 0.272 | 0.264 | 0.162 |
|  | C2 | 484 | 995 | 1358 | 1355 | 903 |
|  |  | 0.095 | 0.195 | 0.267 | 0.266 | 0.177 |

|  |  | Not at all | Minimally | Somewhat | Mostly | Very |
|---|---|---|---|---|---|---|
| A wide range of students respond to the instructor's questions in class. | PC | 674 | 1600 | 2134 | 1552 | 1070 |
|  |  | 0.096 | 0.228 | 0.304 | 0.221 | 0.152 |
|  | C1 | 758 | 1954 | 2104 | 1274 | 709 |
|  |  | 0.111 | 0.287 | 0.309 | 0.187 | 0.104 |
|  | C2 | 384 | 1261 | 1525 | 1141 | 781 |
|  |  | 0.075 | 0.248 | 0.299 | 0.224 | 0.153 |
| The instructor knows my name. | PC | 1225 | 976 | 1212 | 1033 | 2577 |
|  |  | 0.174 | 0.139 | 0.173 | 0.147 | 0.367 |
|  | C1 | 1462 | 981 | 1023 | 896 | 2429 |
|  |  | 0.215 | 0.144 | 0.151 | 0.132 | 0.358 |
|  | C2 | 668 | 521 | 625 | 715 | 2563 |
|  |  | 0.131 | 0.102 | 0.123 | 0.14 | 0.503 |
| *Class is structured to encourage peer-to-peer support among students. | PC | 1223 | 1167 | 1561 | 1446 | 1639 |
|  |  | 0.174 | 0.166 | 0.222 | 0.206 | 0.233 |
|  | C1 | 1354 | 1470 | 1584 | 1306 | 1085 |
|  |  | 0.199 | 0.216 | 0.233 | 0.192 | 0.16 |
|  | C2 | 898 | 965 | 1162 | 985 | 1091 |
|  |  | 0.176 | 0.189 | 0.228 | 0.193 | 0.214 |
| There is a sense of community among the students in my class. | PC | 1110 | 1476 | 1856 | 1405 | 1181 |
|  |  | 0.158 | 0.21 | 0.264 | 0.2 | 0.168 |
|  | C1 | 1033 | 1495 | 1957 | 1380 | 930 |
|  |  | 0.152 | 0.22 | 0.288 | 0.203 | 0.137 |
|  | C2 | 636 | 1010 | 1333 | 1132 | 979 |
|  |  | 0.125 | 0.198 | 0.262 | 0.222 | 0.192 |
| The instructor adjusts teaching based upon what the class understands and does not understand. | PC | 762 | 736 | 1429 | 1810 | 2299 |
|  |  | 0.108 | 0.105 | 0.203 | 0.257 | 0.327 |
|  | C1 | 725 | 961 | 1546 | 1833 | 1730 |
|  |  | 0.107 | 0.141 | 0.228 | 0.27 | 0.255 |
|  | C2 | 359 | 576 | 1043 | 1511 | 1600 |
|  |  | 0.071 | 0.113 | 0.205 | 0.297 | 0.314 |

|  |  | Not at all | Minimally | Somewhat | Mostly | Very |
|---|---|---|---|---|---|---|
| The instructor explains concepts in this class in a variety of ways. | PC | 457 | 844 | 1548 | 1930 | 2256 |
|  |  | 0.065 | 0.12 | 0.22 | 0.274 | 0.321 |
|  | C1 | 420 | 959 | 1583 | 2031 | 1810 |
|  |  | 0.062 | 0.141 | 0.233 | 0.299 | 0.266 |
|  | C2 | 206 | 527 | 1111 | 1621 | 1633 |
|  |  | 0.04 | 0.103 | 0.218 | 0.318 | 0.32 |
| I receive feedback from my instructor on homework, exams, quizzes, etc. | PC | 1006 | 1013 | 1514 | 1539 | 1957 |
|  |  | 0.143 | 0.144 | 0.215 | 0.219 | 0.278 |
|  | C1 | 612 | 1025 | 1544 | 1785 | 1838 |
|  |  | 0.09 | 0.151 | 0.227 | 0.262 | 0.27 |
|  | C2 | 249 | 514 | 953 | 1422 | 1958 |
|  |  | 0.049 | 0.101 | 0.187 | 0.279 | 0.384 |
| I share my ideas (or my group's ideas) during whole class discussions. | PC | 2297 | 1485 | 1455 | 1011 | 782 |
|  |  | 0.327 | 0.211 | 0.207 | 0.144 | 0.111 |
|  | C1 | 2373 | 1638 | 1405 | 843 | 537 |
|  |  | 0.349 | 0.241 | 0.207 | 0.124 | 0.079 |
|  | C2 | 1484 | 1273 | 1073 | 757 | 506 |
|  |  | 0.291 | 0.25 | 0.211 | 0.149 | 0.099 |
| A wide range of students participate in class. | PC | 681 | 1596 | 2205 | 1504 | 1045 |
|  |  | 0.097 | 0.227 | 0.314 | 0.214 | 0.149 |
|  | C1 | 763 | 1986 | 2035 | 1279 | 736 |
|  |  | 0.112 | 0.292 | 0.299 | 0.188 | 0.108 |
|  | C2 | 399 | 1286 | 1485 | 1136 | 794 |
|  |  | 0.078 | 0.252 | 0.291 | 0.223 | 0.156 |
| My instructor uses strategies to encourage participation from a wide range of students. | PC | 755 | 1006 | 1652 | 1785 | 1839 |
|  |  | 0.107 | 0.143 | 0.235 | 0.254 | 0.261 |
|  | C1 | 917 | 1315 | 1791 | 1532 | 1248 |
|  |  | 0.135 | 0.193 | 0.263 | 0.225 | 0.183 |
|  | C2 | 521 | 853 | 1240 | 1198 | 1286 |
|  |  | 0.102 | 0.167 | 0.243 | 0.235 | 0.252 |

## Helpful Course Elements

To what extent are the following course elements helpful to your learning in [course]?

|  |  | Not applicable | Not helpful | Somewhat helpful | Very helpful |
|---|---|---|---|---|---|
| Online homework: | PC | 376<br>0.054 | 906<br>0.129 | 2177<br>0.31 | 3560<br>0.507 |
|  | C1 | 1327<br>0.196 | 771<br>0.114 | 2054<br>0.303 | 2629<br>0.388 |
|  | C2 | 1301<br>0.256 | 536<br>0.106 | 1504<br>0.296 | 1734<br>0.342 |
| Written homework: | PC | 1551<br>0.222 | 822<br>0.118 | 2181<br>0.312 | 2438<br>0.349 |
|  | C1 | 1176<br>0.173 | 572<br>0.084 | 2313<br>0.341 | 2720<br>0.401 |
|  | C2 | 673<br>0.132 | 363<br>0.071 | 1608<br>0.316 | 2442<br>0.48 |
| Exams: | PC | 101<br>0.014 | 1319<br>0.188 | 2969<br>0.424 | 2613<br>0.373 |
|  | C1 | 51<br>0.008 | 1146<br>0.169 | 3027<br>0.446 | 2563<br>0.378 |
|  | C2 | 29<br>0.006 | 724<br>0.142 | 2237<br>0.44 | 2099<br>0.412 |
| Worksheets or handouts in class: | PC | 1270<br>0.181 | 537<br>0.077 | 1895<br>0.27 | 3308<br>0.472 |
|  | C1 | 1508<br>0.222 | 363<br>0.054 | 1929<br>0.284 | 2982<br>0.44 |
|  | C2 | 1055<br>0.207 | 239<br>0.047 | 1279<br>0.251 | 2514<br>0.494 |

## Lab/Recitation Experience

Indicate the degree to which the following statements describe your experience in **recitation/lab sections** of [course]:

|  |  | Not at all | Minimally | Somewhat | Mostly | Very |
|---|---|---|---|---|---|---|
| I listen as the instructor guides me through major topics. | PC | 196 | 124 | 313 | 401 | 515 |
|  |  | 0.127 | 0.08 | 0.202 | 0.259 | 0.332 |
|  | C1 | 234 | 239 | 580 | 993 | 1793 |
|  |  | 0.061 | 0.062 | 0.151 | 0.259 | 0.467 |
|  | C2 | 185 | 126 | 217 | 340 | 618 |
|  |  | 0.124 | 0.085 | 0.146 | 0.229 | 0.416 |
| The lab/recitation activities connect course content to my life and future work. | PC | 329 | 289 | 442 | 266 | 219 |
|  |  | 0.213 | 0.187 | 0.286 | 0.172 | 0.142 |
|  | C1 | 568 | 692 | 1095 | 795 | 688 |
|  |  | 0.148 | 0.18 | 0.285 | 0.207 | 0.179 |
|  | C2 | 322 | 258 | 384 | 286 | 238 |
|  |  | 0.216 | 0.173 | 0.258 | 0.192 | 0.16 |
| *I receive immediate feedback on my work during lab/recitation (e.g., clickers) | PC | 355 | 221 | 314 | 324 | 337 |
|  |  | 0.229 | 0.142 | 0.202 | 0.209 | 0.217 |
|  | C1 | 883 | 489 | 794 | 793 | 885 |
|  |  | 0.23 | 0.127 | 0.207 | 0.206 | 0.23 |
|  | C2 | 444 | 188 | 275 | 282 | 301 |
|  |  | 0.298 | 0.126 | 0.185 | 0.189 | 0.202 |
| I am asked to respond to questions during lab/recitation time. | PC | 437 | 230 | 326 | 280 | 270 |
|  |  | 0.283 | 0.149 | 0.211 | 0.181 | 0.175 |
|  | C1 | 735 | 643 | 835 | 834 | 796 |
|  |  | 0.191 | 0.167 | 0.217 | 0.217 | 0.207 |
|  | C2 | 356 | 231 | 307 | 295 | 303 |
|  |  | 0.239 | 0.155 | 0.206 | 0.198 | 0.203 |
| In my lab/recitation a variety of means (graphs, symbols, tables, etc.) are used to represent course topics and/or solve problems. | PC | 164 | 189 | 381 | 425 | 389 |
|  |  | 0.106 | 0.122 | 0.246 | 0.275 | 0.251 |
|  | C1 | 241 | 420 | 912 | 1105 | 1163 |
|  |  | 0.063 | 0.109 | 0.237 | 0.288 | 0.303 |
|  | C2 | 236 | 201 | 341 | 324 | 389 |
|  |  | 0.158 | 0.135 | 0.229 | 0.217 | 0.261 |

|  |  | Not at all | Minimally | Somewhat | Mostly | Very |
|---|---|---|---|---|---|---|
| I talk with other students about course topics during lab/recitation. | PC | 159 | 155 | 322 | 397 | 516 |
|  |  | 0.103 | 0.1 | 0.208 | 0.256 | 0.333 |
|  | C1 | 410 | 353 | 691 | 955 | 1425 |
|  |  | 0.107 | 0.092 | 0.18 | 0.249 | 0.372 |
|  | C2 | 226 | 128 | 252 | 344 | 536 |
|  |  | 0.152 | 0.086 | 0.17 | 0.231 | 0.361 |
| I constructively criticize other students' ideas during lab/recitation. | PC | 519 | 290 | 328 | 236 | 175 |
|  |  | 0.335 | 0.187 | 0.212 | 0.152 | 0.113 |
|  | C1 | 1122 | 727 | 861 | 641 | 487 |
|  |  | 0.292 | 0.189 | 0.224 | 0.167 | 0.127 |
|  | C2 | 477 | 258 | 312 | 252 | 190 |
|  |  | 0.32 | 0.173 | 0.21 | 0.169 | 0.128 |
| I discuss the difficulties I have with math with other students during lab/recitation. | PC | 174 | 179 | 337 | 392 | 466 |
|  |  | 0.112 | 0.116 | 0.218 | 0.253 | 0.301 |
|  | C1 | 473 | 430 | 767 | 946 | 1228 |
|  |  | 0.123 | 0.112 | 0.2 | 0.246 | 0.319 |
|  | C2 | 259 | 148 | 272 | 355 | 458 |
|  |  | 0.174 | 0.099 | 0.182 | 0.238 | 0.307 |
| I work on problems individually during lab/recitation time. | PC | 153 | 231 | 424 | 349 | 388 |
|  |  | 0.099 | 0.15 | 0.274 | 0.226 | 0.251 |
|  | C1 | 222 | 494 | 985 | 966 | 1171 |
|  |  | 0.058 | 0.129 | 0.257 | 0.252 | 0.305 |
|  | C2 | 203 | 185 | 359 | 359 | 384 |
|  |  | 0.136 | 0.124 | 0.241 | 0.241 | 0.258 |
| I work with other students in small groups during lab/recitation. | PC | 196 | 151 | 242 | 290 | 666 |
|  |  | 0.127 | 0.098 | 0.157 | 0.188 | 0.431 |
|  | C1 | 502 | 375 | 544 | 645 | 1771 |
|  |  | 0.131 | 0.098 | 0.142 | 0.168 | 0.462 |
|  | C2 | 281 | 134 | 156 | 222 | 695 |
|  |  | 0.189 | 0.09 | 0.105 | 0.149 | 0.467 |

|  |  | Not at all | Minimally | Somewhat | Mostly | Very |
|---|---|---|---|---|---|---|
| Multiple approaches to solving a problem are discussed in lab/recitation. | PC | 231 | 206 | 383 | 380 | 344 |
|  |  | 0.15 | 0.133 | 0.248 | 0.246 | 0.223 |
|  | C1 | 309 | 415 | 936 | 1091 | 1083 |
|  |  | 0.081 | 0.108 | 0.244 | 0.285 | 0.282 |
|  | C2 | 230 | 164 | 320 | 388 | 387 |
|  |  | 0.154 | 0.11 | 0.215 | 0.261 | 0.26 |
| I have enough time during lab/recitation to reflect about the processes I use to solve problems. | PC | 150 | 239 | 391 | 430 | 335 |
|  |  | 0.097 | 0.155 | 0.253 | 0.278 | 0.217 |
|  | C1 | 235 | 437 | 923 | 1111 | 1131 |
|  |  | 0.061 | 0.114 | 0.241 | 0.29 | 0.295 |
|  | C2 | 212 | 174 | 343 | 395 | 366 |
|  |  | 0.142 | 0.117 | 0.23 | 0.265 | 0.246 |
| A wide range of students respond to the instructor's questions in lab/recitation. | PC | 333 | 266 | 373 | 308 | 262 |
|  |  | 0.216 | 0.173 | 0.242 | 0.2 | 0.17 |
|  | C1 | 468 | 601 | 994 | 938 | 837 |
|  |  | 0.122 | 0.157 | 0.259 | 0.244 | 0.218 |
|  | C2 | 266 | 241 | 358 | 335 | 291 |
|  |  | 0.178 | 0.162 | 0.24 | 0.225 | 0.195 |
| The instructor knows my name. | PC | 521 | 222 | 251 | 181 | 369 |
|  |  | 0.337 | 0.144 | 0.163 | 0.117 | 0.239 |
|  | C1 | 644 | 406 | 571 | 630 | 1587 |
|  |  | 0.168 | 0.106 | 0.149 | 0.164 | 0.413 |
|  | C2 | 271 | 179 | 188 | 194 | 660 |
|  |  | 0.182 | 0.12 | 0.126 | 0.13 | 0.442 |
| Lab/recitation is structured to encourage peer-to-peer support among students. | PC | 151 | 171 | 293 | 386 | 546 |
|  |  | 0.098 | 0.111 | 0.189 | 0.25 | 0.353 |
|  | C1 | 355 | 352 | 696 | 901 | 1541 |
|  |  | 0.092 | 0.092 | 0.181 | 0.234 | 0.401 |
|  | C2 | 229 | 134 | 245 | 314 | 569 |
|  |  | 0.154 | 0.09 | 0.164 | 0.211 | 0.382 |

|  |  | Not at all | Minimally | Somewhat | Mostly | Very |
|---|---|---|---|---|---|---|
| There is a sense of community among the students in my lab/recitation. | PC | 258 | 247 | 393 | 340 | 309 |
|  |  | 0.167 | 0.16 | 0.254 | 0.22 | 0.2 |
|  | C1 | 387 | 487 | 958 | 948 | 1056 |
|  |  | 0.101 | 0.127 | 0.25 | 0.247 | 0.275 |
|  | C2 | 240 | 165 | 351 | 361 | 373 |
|  |  | 0.161 | 0.111 | 0.236 | 0.242 | 0.25 |
| The instructor adjusts teaching based on what the class understands and does not understand. | PC | 289 | 209 | 335 | 336 | 375 |
|  |  | 0.187 | 0.135 | 0.217 | 0.218 | 0.243 |
|  | C1 | 385 | 364 | 706 | 1052 | 1333 |
|  |  | 0.1 | 0.095 | 0.184 | 0.274 | 0.347 |
|  | C2 | 253 | 149 | 279 | 370 | 440 |
|  |  | 0.17 | 0.1 | 0.187 | 0.248 | 0.295 |
| The instructor explains concepts in this lab/recitation in a variety of ways. | PC | 249 | 217 | 368 | 363 | 348 |
|  |  | 0.161 | 0.14 | 0.238 | 0.235 | 0.225 |
|  | C1 | 313 | 424 | 866 | 1062 | 1175 |
|  |  | 0.082 | 0.11 | 0.226 | 0.277 | 0.306 |
|  | C2 | 233 | 190 | 324 | 376 | 366 |
|  |  | 0.156 | 0.128 | 0.218 | 0.253 | 0.246 |
| I receive feedback from my instructor on homework, exams, quizzes, etc. | PC | 260 | 210 | 306 | 376 | 397 |
|  |  | 0.168 | 0.136 | 0.198 | 0.243 | 0.256 |
|  | C1 | 322 | 407 | 774 | 1045 | 1293 |
|  |  | 0.084 | 0.106 | 0.202 | 0.272 | 0.337 |
|  | C2 | 217 | 143 | 261 | 338 | 530 |
|  |  | 0.146 | 0.096 | 0.175 | 0.227 | 0.356 |
| I share my ideas (or my group's ideas) during whole lab/recitation discussions. | PC | 457 | 251 | 286 | 290 | 261 |
|  |  | 0.296 | 0.162 | 0.185 | 0.188 | 0.169 |
|  | C1 | 849 | 653 | 858 | 722 | 756 |
|  |  | 0.221 | 0.17 | 0.224 | 0.188 | 0.197 |
|  | C2 | 396 | 228 | 322 | 276 | 267 |
|  |  | 0.266 | 0.153 | 0.216 | 0.185 | 0.179 |

|  |  | Not at all | Minimally | Somewhat | Mostly | Very |
|---|---|---|---|---|---|---|
| A wide range of students participate in lab/recitation. | PC | 200 | 269 | 358 | 357 | 358 |
|  |  | 0.13 | 0.174 | 0.232 | 0.232 | 0.232 |
|  | C1 | 313 | 551 | 938 | 994 | 1041 |
|  |  | 0.082 | 0.144 | 0.244 | 0.259 | 0.271 |
|  | C2 | 213 | 195 | 354 | 371 | 358 |
|  |  | 0.143 | 0.131 | 0.237 | 0.249 | 0.24 |
| My instructor uses strategies to encourage participation from a wide range of students. | PC | 265 | 220 | 392 | 337 | 334 |
|  |  | 0.171 | 0.142 | 0.253 | 0.218 | 0.216 |
|  | C1 | 448 | 499 | 922 | 945 | 1017 |
|  |  | 0.117 | 0.13 | 0.241 | 0.247 | 0.265 |
|  | C2 | 282 | 205 | 337 | 299 | 369 |
|  |  | 0.189 | 0.137 | 0.226 | 0.2 | 0.247 |

## Class Comparisons

Consider your regular course meetings and primary instructor of [course]. As compared to other students in class…

|  |  | A lot less | Somewhat less | The same | Somewhat more | A lot more |
|---|---|---|---|---|---|---|
| How much opportunity do you get to answer questions in class? | PC | 516 | 478 | 5366 | 371 | 218 |
|  |  | 0.074 | 0.069 | 0.772 | 0.053 | 0.031 |
|  | C1 | 459 | 428 | 5370 | 331 | 128 |
|  |  | 0.068 | 0.064 | 0.8 | 0.049 | 0.019 |
|  | C2 | 252 | 297 | 4118 | 303 | 90 |
|  |  | 0.05 | 0.059 | 0.814 | 0.06 | 0.018 |
| How much attention does the instructor give to your questions? | PC | 346 | 330 | 5812 | 301 | 148 |
|  |  | 0.05 | 0.048 | 0.838 | 0.043 | 0.021 |
|  | C1 | 279 | 316 | 5793 | 225 | 92 |
|  |  | 0.042 | 0.047 | 0.864 | 0.034 | 0.014 |
|  | C2 | 145 | 213 | 4417 | 214 | 68 |
|  |  | 0.029 | 0.042 | 0.873 | 0.042 | 0.013 |

|  |  | A lot less | Somewhat less | The same | Somewhat more | A lot more |
|---|---|---|---|---|---|---|
| How much help do you get from the instructor? | PC | 492 | 648 | 5133 | 462 | 204 |
|  |  | 0.071 | 0.093 | 0.74 | 0.067 | 0.029 |
|  | C1 | 446 | 664 | 5024 | 471 | 92 |
|  |  | 0.067 | 0.099 | 0.75 | 0.07 | 0.014 |
|  | C2 | 235 | 466 | 3823 | 438 | 93 |
|  |  | 0.046 | 0.092 | 0.756 | 0.087 | 0.018 |
| How much encouragement do you receive from the instructor? | PC | 416 | 447 | 5402 | 460 | 207 |
|  |  | 0.06 | 0.064 | 0.779 | 0.066 | 0.03 |
|  | C1 | 382 | 388 | 5424 | 402 | 103 |
|  |  | 0.057 | 0.058 | 0.81 | 0.06 | 0.015 |
|  | C2 | 168 | 264 | 4164 | 351 | 103 |
|  |  | 0.033 | 0.052 | 0.825 | 0.07 | 0.02 |
| How much opportunity do you get to contribute to class discussions? | PC | 409 | 401 | 5575 | 359 | 189 |
|  |  | 0.059 | 0.058 | 0.804 | 0.052 | 0.027 |
|  | C1 | 360 | 382 | 5545 | 320 | 99 |
|  |  | 0.054 | 0.057 | 0.827 | 0.048 | 0.015 |
|  | C2 | 181 | 295 | 4216 | 274 | 86 |
|  |  | 0.036 | 0.058 | 0.835 | 0.054 | 0.017 |
| How much praise does your work receive? | PC | 513 | 445 | 5467 | 354 | 144 |
|  |  | 0.074 | 0.064 | 0.79 | 0.051 | 0.021 |
|  | C1 | 446 | 457 | 5414 | 306 | 72 |
|  |  | 0.067 | 0.068 | 0.809 | 0.046 | 0.011 |
|  | C2 | 217 | 327 | 4174 | 263 | 65 |
|  |  | 0.043 | 0.065 | 0.827 | 0.052 | 0.013 |

## Lab/Recitation Comparisons

Consider your **recitation/lab section and recitation/lab instructor**. As compared to other students in class…

|  |  | A lot less | Somewhat less | The same | Somewhat more | A lot more |
|---|---|---|---|---|---|---|
| How much opportunity do you get to answer questions in lab/recitation? | PC | 85<br>0.055 | 90<br>0.058 | 1256<br>0.816 | 80<br>0.052 | 28<br>0.018 |
|  | C1 | 138<br>0.036 | 164<br>0.043 | 3204<br>0.836 | 281<br>0.073 | 45<br>0.012 |
|  | C2 | 50<br>0.034 | 54<br>0.036 | 1272<br>0.857 | 84<br>0.057 | 25<br>0.017 |
| How much attention does the instructor give your questions? | PC | 69<br>0.045 | 93<br>0.06 | 1259<br>0.818 | 91<br>0.059 | 27<br>0.018 |
|  | C1 | 99<br>0.026 | 162<br>0.042 | 3251<br>0.849 | 272<br>0.071 | 45<br>0.012 |
|  | C2 | 50<br>0.034 | 50<br>0.034 | 1274<br>0.86 | 82<br>0.055 | 26<br>0.018 |
| How much help do you get from the instructor? | PC | 81<br>0.053 | 118<br>0.077 | 1174<br>0.763 | 123<br>0.08 | 42<br>0.027 |
|  | C1 | 141<br>0.037 | 247<br>0.065 | 2983<br>0.779 | 369<br>0.096 | 87<br>0.023 |
|  | C2 | 62<br>0.042 | 88<br>0.059 | 1180<br>0.796 | 121<br>0.082 | 32<br>0.022 |
| How much encouragement do you receive from the instructor? | PC | 75<br>0.049 | 99<br>0.065 | 1214<br>0.791 | 109<br>0.071 | 37<br>0.024 |
|  | C1 | 112<br>0.029 | 166<br>0.043 | 3194<br>0.836 | 280<br>0.073 | 69<br>0.018 |
|  | C2 | 53<br>0.036 | 61<br>0.041 | 1253<br>0.847 | 83<br>0.056 | 30<br>0.02 |
| How much opportunity do you get to contribute to lab/recitation discussion? | PC | 74<br>0.048 | 91<br>0.059 | 1268<br>0.827 | 75<br>0.049 | 26<br>0.017 |
|  | C1 | 129<br>0.034 | 149<br>0.039 | 3295<br>0.861 | 195<br>0.051 | 57<br>0.015 |
|  | C2 | 52<br>0.035 | 49<br>0.033 | 1288<br>0.87 | 70<br>0.047 | 22<br>0.015 |

|  |  | A lot less | Somewhat less | The same | Somewhat more | A lot more |
|---|---|---|---|---|---|---|
| How much praise does your work receive? | PC | 97 | 84 | 1227 | 83 | 29 |
|  |  | 0.064 | 0.055 | 0.807 | 0.055 | 0.019 |
|  | C1 | 144 | 179 | 3215 | 215 | 54 |
|  |  | 0.038 | 0.047 | 0.844 | 0.056 | 0.014 |
|  | C2 | 56 | 69 | 1256 | 69 | 25 |
|  |  | 0.038 | 0.047 | 0.852 | 0.047 | 0.017 |

## Classroom Climate

How would you describe the overall climate in [course]?

|  | (1) Excluding and Hostile | (2) | (3) | (4) | (5) Including and Friendly |
|---|---|---|---|---|---|
| PC | 171 | 337 | 1140 | 1759 | 3509 |
|  | 0.025 | 0.049 | 0.165 | 0.254 | 0.507 |
| C1 | 139 | 348 | 1223 | 2014 | 2979 |
|  | 0.021 | 0.052 | 0.182 | 0.3 | 0.444 |
| C2 | 90 | 223 | 791 | 1379 | 2567 |
|  | 0.018 | 0.044 | 0.157 | 0.273 | 0.508 |
|  | (1) Intellectually Boring | (2) | (3) | (4) | (5) Intellectually Engaging |
| PC | 631 | 903 | 1743 | 1779 | 1860 |
|  | 0.091 | 0.131 | 0.252 | 0.257 | 0.269 |
| C1 | 536 | 886 | 1598 | 1951 | 1731 |
|  | 0.08 | 0.132 | 0.238 | 0.291 | 0.258 |
| C2 | 236 | 477 | 1044 | 1580 | 1711 |
|  | 0.047 | 0.094 | 0.207 | 0.313 | 0.339 |
|  | (1) Academically Easy | (2) | (3) | (4) | (5) Academically Rigorous |
| PC | 272 | 717 | 2190 | 2017 | 1711 |
|  | 0.039 | 0.104 | 0.317 | 0.292 | 0.248 |
| C1 | 192 | 457 | 1592 | 2211 | 2253 |
|  | 0.029 | 0.068 | 0.237 | 0.33 | 0.336 |
| C2 | 72 | 180 | 835 | 1737 | 2222 |
|  | 0.014 | 0.036 | 0.165 | 0.344 | 0.44 |

## Malleability of Math Ability

I believe that my math ability can be improved through dedication and hard work.

*Response options: Strongly agree (1) to Strongly disagree (6).*

|  |  | 1 (SA) | 2 | 3 | 4 | 5 | 6 (SD) |
|---|---|---|---|---|---|---|---|
| I believe that my math ability can be improved through dedication and hard work | PC | 3993 | 2002 | 612 | 133 | 44 | 129 |
|  |  | 0.578 | 0.29 | 0.089 | 0.019 | 0.006 | 0.019 |
|  | C1 | 3711 | 2202 | 520 | 103 | 46 | 97 |
|  |  | 0.556 | 0.33 | 0.078 | 0.015 | 0.007 | 0.015 |
|  | C2 | 3026 | 1482 | 366 | 68 | 24 | 68 |
|  |  | 0.601 | 0.294 | 0.073 | 0.014 | 0.005 | 0.014 |

## Mathematical Attitudes

Please indicate your level of agreement for the following statements from the beginning of the course and now.

*Response options: Strongly agree (1) to Strongly disagree (6).*

| **I am interested in mathematics:** | | | | | | | | | | | | |
|---|---|---|---|---|---|---|---|---|---|---|---|---|
| | Attitude at the beginning of the course | | | | | | Current attitude in the course | | | | | |
| | 1 (SA) | 2 | 3 | 4 | 5 | 6 (SD) | 1 (SA) | 2 | 3 | 4 | 5 | 6 (SD) |
| PC | 1352 | 1752 | 1396 | 678 | 972 | 797 | 1240 | 1722 | 1437 | 733 | 768 | 821 |
| | 0.195 | 0.252 | 0.201 | 0.098 | 0.14 | 0.115 | 0.184 | 0.256 | 0.214 | 0.109 | 0.114 | 0.122 |
| C1 | 1640 | 2117 | 1472 | 506 | 607 | 372 | 1378 | 2018 | 1434 | 618 | 575 | 527 |
| | 0.244 | 0.315 | 0.219 | 0.075 | 0.09 | 0.055 | 0.21 | 0.308 | 0.219 | 0.094 | 0.088 | 0.08 |
| C2 | 1507 | 1840 | 1005 | 309 | 265 | 144 | 1359 | 1693 | 1017 | 381 | 299 | 228 |
| | 0.297 | 0.363 | 0.198 | 0.061 | 0.052 | 0.028 | 0.273 | 0.34 | 0.204 | 0.077 | 0.06 | 0.046 |
| **I enjoy doing mathematics:** | | | | | | | | | | | | |
| | Attitude at the beginning of the course | | | | | | Current attitude in the course | | | | | |
| | 1 (SA) | 2 | 3 | 4 | 5 | 6 (SD) | 1 (SA) | 2 | 3 | 4 | 5 | 6 (SD) |
| PC | 1161 | 1639 | 1365 | 773 | 970 | 1029 | 1061 | 1606 | 1470 | 808 | 806 | 949 |
| | 0.167 | 0.236 | 0.197 | 0.111 | 0.14 | 0.148 | 0.158 | 0.24 | 0.219 | 0.121 | 0.12 | 0.142 |
| C1 | 1452 | 1987 | 1504 | 696 | 597 | 467 | 1154 | 1815 | 1554 | 782 | 638 | 587 |
| | 0.217 | 0.296 | 0.224 | 0.104 | 0.089 | 0.07 | 0.177 | 0.278 | 0.238 | 0.12 | 0.098 | 0.09 |
| C2 | 1357 | 1692 | 1068 | 446 | 306 | 196 | 1103 | 1515 | 1172 | 514 | 364 | 304 |
| | 0.268 | 0.334 | 0.211 | 0.088 | 0.06 | 0.039 | 0.222 | 0.305 | 0.236 | 0.103 | 0.073 | 0.061 |

| | \multicolumn{6}{l|}{**I am confident in my mathematical abilities:**} | | | | | | |
|---|---|---|---|---|---|---|---|---|---|---|---|---|
| | Attitude at the beginning of the course | | | | | | Current attitude in the course | | | | | |
| | 1 (SA) | 2 | 3 | 4 | 5 | 6 (SD) | 1 (SA) | 2 | 3 | 4 | 5 | 6 (SD) |
| PC | 1068 | 1859 | 1717 | 845 | 755 | 668 | 1090 | 2009 | 1654 | 755 | 541 | 553 |
|    | 0.155 | 0.269 | 0.248 | 0.122 | 0.109 | 0.097 | 0.165 | 0.304 | 0.251 | 0.114 | 0.082 | 0.084 |
| C1 | 1408 | 2134 | 1561 | 724 | 491 | 354 | 1062 | 1883 | 1654 | 802 | 566 | 450 |
|    | 0.211 | 0.32 | 0.234 | 0.109 | 0.074 | 0.053 | 0.165 | 0.293 | 0.258 | 0.125 | 0.088 | 0.07 |
| C2 | 1216 | 1801 | 1150 | 451 | 264 | 166 | 858 | 1586 | 1289 | 510 | 385 | 277 |
|    | 0.241 | 0.357 | 0.228 | 0.089 | 0.052 | 0.033 | 0.175 | 0.323 | 0.263 | 0.104 | 0.078 | 0.056 |
| | \multicolumn{6}{l|}{**I am able to learn mathematics:**} | | | | | | |
| | Attitude at the beginning of the course | | | | | | Current attitude in the course | | | | | |
| | 1 (SA) | 2 | 3 | 4 | 5 | 6 (SD) | 1 (SA) | 2 | 3 | 4 | 5 | 6 (SD) |
| PC | 2004 | 2654 | 1165 | 386 | 280 | 217 | 2292 | 2527 | 951 | 314 | 209 | 238 |
|    | 0.299 | 0.396 | 0.174 | 0.058 | 0.042 | 0.032 | 0.351 | 0.387 | 0.146 | 0.048 | 0.032 | 0.036 |
| C1 | 2167 | 2769 | 1039 | 271 | 137 | 99 | 2112 | 2499 | 1069 | 341 | 176 | 163 |
|    | 0.334 | 0.427 | 0.16 | 0.042 | 0.021 | 0.015 | 0.332 | 0.393 | 0.168 | 0.054 | 0.028 | 0.026 |
| C2 | 1918 | 2073 | 671 | 159 | 72 | 43 | 1848 | 1866 | 724 | 224 | 114 | 88 |
|    | 0.389 | 0.42 | 0.136 | 0.032 | 0.015 | 0.009 | 0.38 | 0.384 | 0.149 | 0.046 | 0.023 | 0.018 |

# U-GPIPS: Student Instructor Survey

The U-GPIPS data contain 331 responses from undergraduate and graduate student instructors from twelve universities across the United States. Of these responses, 294 come from graduate student instructors, 34 from undergraduate student instructors, and 3 from student instructors who did not report their current degree program. We first present the self-identified demographic information of these students.

As noted in the introduction to this report, instructors responded for each course they were teaching in a particular term. In the identity and context portion of this report, instructors are counted only once per term; in the course-specific information their responses are reported separately for each *course* they were teaching at the time they completed the survey.

## Student Instructor Identity & Individual Context

These questions asked about student instructors' identity, including demographic information and academic context. Questions about Gender, Race and Ethnicity, and Sexual Orientation were presented as "select all that apply" items with the option to self-identify descriptors not on our lists. We intentionally provided more nuanced options than are generally use in demographic reports (e.g., the US census) as we value the complex nature of identity. Responses to these items were recoded for interpretability in this report; details are provided alongside each item.

### Degree and Career Aspirations

What degree(s) or certifications do you intend to obtain from [your institution]? Mark all that apply.

| | |
|---|---|
| B.A./B.S. | 50 (0.151) |
| Teaching certification/credential | 8 (0.024) |
| M.A./M.S. in mathematics | 81 (0.245) |
| M.A./M.S. in mathematics education | 6 (0.018) |
| Ph.D. in mathematics | 203 (0.613) |
| Ph.D. in mathematics education | 30 (0.091) |
| Ed.D. | 1 (0.003) |
| Other | 11 (0.033) |

What is your intended career trajectory?

| | |
|---|---|
| Academic position at a 4-year college or university – teaching focused | 93 (0.347) |
| Academic position at a 4-year college or university – research focused | 112 (0.418) |
| Academic position at a 2-year college | 5 (0.019) |

| Non-academic position (industry, government, etc.) | 34 (0.127) |
|---|---|
| Other (please explain) | 24 (0.09) |

## Gender

Do you consider yourself to be (Select all that apply):

| Cisgender Man | Cisgender Woman | Gender Non-binary | Transgender | Transgender Man | Transgender Woman | Not reported |
|---|---|---|---|---|---|---|
| 165 | 113 | 1 | 2 | 0 | 0 | 50 |
| 0.498 | 0.341 | 0.003 | 0.006 | 0.000 | 0.000 | 0.151 |

*Survey response options: Man, Transgender, Woman, Not listed (please specify), Prefer not to disclose.*

*Selections of only Man, only Woman, or only Transgender are coded as Cisgender Man, Cisgender Woman, and Transgender respectively; selections of Transgender and Man or Transgender and Woman are coded as Transgender Man and Transgender Woman, respectively; selections of Man and Woman (with or without additional selections) or self-described non-binary gender identities (e.g., genderfluid, agender) are coded as Gender Non-Binary; selections of Prefer not to disclose and no selection were coded as Not reported.*

## Race and Ethnicity

Do you consider yourself to be (Select all that apply):

| Asian | 62 (0.187) |
|---|---|
| Alaska Native or Native American | 0 (0.000) |
| Black or African American | 8 (0.024) |
| Hispanic or Latinx | 8 (0.024) |
| Middle Eastern or North African | 5 (0.015) |
| Multiple Race | 16 (0.048) |
| Native Hawaiian or Pacific Islander | 0 (0.000) |
| White | 180 (0.544) |
| Not reported | 52 (0.157) |

*Survey response options: Alaskan Native or Native American, Black or African American, Central Asian, East Asian, Hispanic or Latinx, Middle Eastern or North African, Native Hawaiian or Pacific Islander, Southeast Asian, South Asian, White, Not listed (please specify), Prefer not to disclose.*

*Selections of more than one of Central Asian, East Asian, Southeast Asian, South Asian without any other choices are coded as Asian; selections of multiple responses (aside from*

*the previous) are coded as Multiple Race/Ethnicity; selections of Prefer not to disclose and no selection are coded as Not reported.*

### Sexual Orientation

Do you consider yourself to be (Select all that apply):

| Straight | Gay | Lesbian | Bisexual | Asexual | Queer | Queer+ | Straight+ | Not reported |
|---|---|---|---|---|---|---|---|---|
| 205 | 7 | 0 | 19 | 5 | 1 | 3 | 2 | 89 |
| 0.619 | 0.021 | 0.000 | 0.057 | 0.015 | 0.003 | 0.009 | 0.006 | 0.269 |

*Survey response options: Asexual, Bisexual, Gay, Lesbian, Queer, Straight (Heterosexual), Not listed (please specify), Prefer not to disclose.*

*Categories are based on research from Voigt (2020) where selections of exactly one listed option are not recoded; selections of Straight (Heterosexual) and one or more additional options are coded as Straight+; selections of more than one option but not Straight (Heterosexual) are coded as Queer+; selections of Prefer not to disclose and no selection are coded as Not reported.*

### Additional Context / Identity Markers

Do you consider yourself to be (Select all that apply):

| | |
|---|---|
| International Student | 80 (0.242) |
| First-generation college student | 72 (0.218) |
| First-generation higher education | 87 (0.263) |
| Student with a disability | 6 (0.018) |
| English language learner | 45 (0.136) |
| Parent or guardian | 14 (0.042) |

## Course Specific Information

This section describes the participants' perceptions of their role in Precalculus, Calculus I, and/or Calculus II (it may be the case that one participant taught more than one of the above courses during the semester). Each result will be shown in the order Precalculus, Calculus I, Calculus II.

### Course Role

What is your official role with regards to [course]?

|    | Instructor of record/ primary instructor | Recitation/discussion/lab section leader | Other |
|----|---|---|---|
| PC | 44      | 12    | 3     |
|    | 0.746   | 0.203 | 0.051 |
| C1 | 77      | 82    | 1     |
|    | 0.481   | 0.512 | 0.006 |
| C2 | 56      | 37    | 0     |
|    | 0.602   | 0.398 | 0.000 |

### Decisions about the course

How are most decisions about course content (e.g. syllabi, exams, homework, pacing, grading) made for [course]? Clarify if you wish.

|    | I make most decisions | I am part of a team that makes most decisions | Someone else makes most decisions |
|----|---|---|---|
| PC | 1       | 25    | 32    |
|    | 0.017   | 0.431 | 0.552 |
| C1 | 17      | 40    | 103   |
|    | 0.106   | 0.25  | 0.644 |
| C2 | 21      | 22    | 50    |
|    | 0.226   | 0.237 | 0.538 |

How are most decisions about instructional approach (e.g., use of clickers, group work, active learning) made for [course]? Clarify if you wish.

|    | I make most decisions | I am part of a team that makes most decisions | Someone else makes most decisions |
|----|---|---|---|
| PC | 21      | 16    | 22    |
|    | 0.356   | 0.271 | 0.373 |
| C1 | 71      | 37    | 50    |
|    | 0.449   | 0.234 | 0.316 |
| C2 | 55      | 15    | 24    |
|    | 0.585   | 0.16  | 0.255 |

## Percentages of Class Time

When you are teaching [course], what percent of that time do students spend… [must total 100]

|  | n | mean | sd | median |
|---|---|---|---|---|
| **All Student Instructors** | | | | |
| *Working on tasks individually* | 314 | 10.98 | 12.67 | 10 |
| *Working on tasks in small groups* | 314 | 29.81 | 27.82 | 20 |
| *Engaging in whole class discussion* | 314 | 12.46 | 13.49 | 10 |
| *Listening to the instructor lecture or solve problems* | 314 | 46.75 | 27.79 | 50 |
| **Precalculus Instructors** | | | | |
| *Working on tasks individually* | 59 | 13.2 | 18.77 | 10 |
| *Working on tasks in small groups* | 59 | 25.39 | 21.89 | 23 |
| *Engaging in whole class discussion* | 59 | 18.31 | 17.05 | 15 |
| *Listening to the instructor lecture or solve problems* | 59 | 43.1 | 22.56 | 41 |
| **Calculus 1 Instructors** | | | | |
| *Working on tasks individually* | 161 | 11.3 | 11.1 | 10 |
| *Working on tasks in small groups* | 161 | 31.29 | 27.98 | 24 |
| *Engaging in whole class discussion* | 161 | 12.4 | 13.51 | 10 |
| *Listening to the instructor lecture or solve problems* | 161 | 45 | 28.16 | 45 |
| **Calculus 2 Instructors** | | | | |
| *Working on tasks individually* | 94 | 9.043 | 10.07 | 6 |
| *Working on tasks in small groups* | 94 | 30.03 | 30.71 | 20 |
| *Engaging in whole class discussion* | 94 | 8.883 | 9.065 | 6 |
| *Listening to the instructor lecture or solve problems* | 94 | 52.04 | 29.58 | 60 |

## Classroom Experience

Please indicate the degree to which the following statements are descriptive of your teaching.

| | | Not at all | Minimally | Somewhat | Mostly | Very |
|---|---|---|---|---|---|---|
| I guide students through major topics as they listen: | PC | 3 | 7 | 8 | 27 | 12 |
| | | 0.053 | 0.123 | 0.14 | 0.474 | 0.211 |
| | C1 | 1 | 9 | 32 | 58 | 52 |
| | | 0.007 | 0.059 | 0.211 | 0.382 | 0.342 |
| | C2 | 1 | 7 | 15 | 36 | 31 |
| | | 0.011 | 0.078 | 0.167 | 0.4 | 0.344 |

|  |  | Not at all | Minimally | Somewhat | Mostly | Very |
|---|---|---|---|---|---|---|
| I provide activities that connect course content to my students' lives and future work: | PC | 6 | 16 | 18 | 7 | 10 |
| | | 0.105 | 0.281 | 0.316 | 0.123 | 0.175 |
| | C1 | 24 | 37 | 48 | 31 | 15 |
| | | 0.155 | 0.239 | 0.31 | 0.2 | 0.097 |
| | C2 | 11 | 21 | 32 | 19 | 8 |
| | | 0.121 | 0.231 | 0.352 | 0.209 | 0.088 |
| My syllabus contains the specific topics that will be covered in every class session: | PC | 8 | 9 | 7 | 9 | 23 |
| | | 0.143 | 0.161 | 0.125 | 0.161 | 0.411 |
| | C1 | 18 | 16 | 17 | 37 | 60 |
| | | 0.122 | 0.108 | 0.115 | 0.25 | 0.405 |
| | C2 | 12 | 9 | 13 | 21 | 36 |
| | | 0.132 | 0.099 | 0.143 | 0.231 | 0.396 |
| I provide students with immediate feedback on their word during class (e.g., student response systems; short quizzes): | PC | 14 | 4 | 13 | 20 | 6 |
| | | 0.246 | 0.07 | 0.228 | 0.351 | 0.105 |
| | C1 | 23 | 29 | 30 | 34 | 38 |
| | | 0.149 | 0.188 | 0.195 | 0.221 | 0.247 |
| | C2 | 15 | 22 | 11 | 23 | 20 |
| | | 0.165 | 0.242 | 0.121 | 0.253 | 0.22 |
| I structure my course with the assumption that most of the students have little useful knowledge of the topics: | PC | 2 | 7 | 16 | 17 | 14 |
| | | 0.036 | 0.125 | 0.286 | 0.304 | 0.25 |
| | C1 | 2 | 17 | 43 | 56 | 35 |
| | | 0.013 | 0.111 | 0.281 | 0.366 | 0.229 |
| | C2 | 2 | 13 | 21 | 31 | 22 |
| | | 0.022 | 0.146 | 0.236 | 0.348 | 0.247 |
| I use student assessment results to guide the direction of my instruction during the semester: | PC | 6 | 8 | 12 | 16 | 15 |
| | | 0.105 | 0.14 | 0.211 | 0.281 | 0.263 |
| | C1 | 16 | 18 | 39 | 41 | 38 |
| | | 0.105 | 0.118 | 0.257 | 0.27 | 0.25 |
| | C2 | 12 | 13 | 20 | 24 | 22 |
| | | 0.132 | 0.143 | 0.22 | 0.264 | 0.242 |

|  |  | Not at all | Minimally | Somewhat | Mostly | Very |
|---|---|---|---|---|---|---|
| I ask students to respond to questions during class time: | PC | 1 | 5 | 13 | 38 | 1 |
| | | 0.018 | 0.088 | 0.228 | 0.667 | 0.018 |
| | C1 | 5 | 21 | 42 | 89 | 5 |
| | | 0.032 | 0.134 | 0.268 | 0.567 | 0.032 |
| | C2 | 1 | 5 | 10 | 25 | 50 |
| | | 0.011 | 0.055 | 0.11 | 0.275 | 0.549 |
| I use student questions and comments to determine the focus and direction of classroom lessons: | PC | 2 | 7 | 13 | 21 | 14 |
| | | 0.035 | 0.123 | 0.228 | 0.368 | 0.246 |
| | C1 | 6 | 13 | 37 | 56 | 42 |
| | | 0.039 | 0.084 | 0.24 | 0.364 | 0.273 |
| | C2 | 2 | 5 | 27 | 36 | 21 |
| | | 0.022 | 0.055 | 0.297 | 0.396 | 0.231 |
| In my class a variety of means (models, drawings, graphs, symbols, simulations, tables, etc.) are used to represent course topics and/or solve problems: | PC | 1 | 3 | 12 | 23 | 18 |
| | | 0.018 | 0.053 | 0.211 | 0.404 | 0.316 |
| | C1 | 4 | 13 | 35 | 47 | 58 |
| | | 0.025 | 0.083 | 0.223 | 0.299 | 0.369 |
| | C2 | 6 | 3 | 27 | 30 | 25 |
| | | 0.066 | 0.033 | 0.297 | 0.33 | 0.275 |
| I structure class so that students explore or discuss their understanding of concepts before direct instruction: | PC | 8 | 22 | 13 | 11 | 3 |
| | | 0.14 | 0.386 | 0.228 | 0.193 | 0.053 |
| | C1 | 24 | 35 | 45 | 31 | 20 |
| | | 0.155 | 0.226 | 0.29 | 0.2 | 0.129 |
| | C2 | 20 | 18 | 25 | 21 | 7 |
| | | 0.22 | 0.198 | 0.275 | 0.231 | 0.077 |
| My class sessions are structured to give students a clear/structured set of notes: | PC | 4 | 5 | 7 | 21 | 20 |
| | | 0.07 | 0.088 | 0.123 | 0.368 | 0.351 |
| | C1 | 16 | 15 | 31 | 41 | 50 |
| | | 0.105 | 0.098 | 0.203 | 0.268 | 0.327 |
| | C2 | 5 | 9 | 16 | 34 | 27 |
| | | 0.055 | 0.099 | 0.176 | 0.374 | 0.297 |

|  |  | Not at all | Minimally | Somewhat | Mostly | Very |
|---|---|---|---|---|---|---|
| I structure class so that students talk with one another about course topics: | PC | 2 | 10 | 12 | 18 | 14 |
|  |  | 0.036 | 0.179 | 0.214 | 0.321 | 0.25 |
|  | C1 | 11 | 21 | 36 | 37 | 52 |
|  |  | 0.07 | 0.134 | 0.229 | 0.236 | 0.331 |
|  | C2 | 7 | 9 | 26 | 30 | 19 |
|  |  | 0.077 | 0.099 | 0.286 | 0.33 | 0.209 |
| I structure class so that students constructively criticize one another's ideas: | PC | 10 | 13 | 23 | 8 | 3 |
|  |  | 0.175 | 0.228 | 0.404 | 0.14 | 0.053 |
|  | C1 | 42 | 41 | 43 | 22 | 9 |
|  |  | 0.268 | 0.261 | 0.274 | 0.14 | 0.057 |
|  | C2 | 19 | 25 | 24 | 17 | 6 |
|  |  | 0.209 | 0.275 | 0.264 | 0.187 | 0.066 |
| I structure class so that students discuss their mathematical difficulties with other students: | PC | 7 | 7 | 17 | 19 | 8 |
|  |  | 0.121 | 0.121 | 0.293 | 0.328 | 0.138 |
|  | C1 | 20 | 21 | 48 | 33 | 34 |
|  |  | 0.128 | 0.135 | 0.308 | 0.212 | 0.218 |
|  | C2 | 10 | 15 | 27 | 27 | 12 |
|  |  | 0.11 | 0.165 | 0.297 | 0.297 | 0.132 |
| I structure class so that students work on problems individually during class: | PC | 13 | 9 | 20 | 11 | 4 |
|  |  | 0.228 | 0.158 | 0.351 | 0.193 | 0.07 |
|  | C1 | 29 | 45 | 48 | 19 | 15 |
|  |  | 0.186 | 0.288 | 0.308 | 0.122 | 0.096 |
|  | C2 | 18 | 27 | 17 | 18 | 11 |
|  |  | 0.198 | 0.297 | 0.187 | 0.198 | 0.121 |
| I structure class so that students work together in pairs or small groups: | PC | 4 | 7 | 11 | 16 | 20 |
|  |  | 0.069 | 0.121 | 0.19 | 0.276 | 0.345 |
|  | C1 | 13 | 18 | 30 | 34 | 62 |
|  |  | 0.083 | 0.115 | 0.191 | 0.217 | 0.395 |
|  | C2 | 9 | 15 | 17 | 18 | 32 |
|  |  | 0.099 | 0.165 | 0.187 | 0.198 | 0.352 |

|  |  | Not at all | Minimally | Somewhat | Mostly | Very |
|---|---|---|---|---|---|---|
| I structure class so that more than one approach to solving a problem is discussed: | PC | 2 | 2 | 11 | 25 | 17 |
|  |  | 0.035 | 0.035 | 0.193 | 0.439 | 0.298 |
|  | C1 | 1 | 12 | 42 | 57 | 43 |
|  |  | 0.006 | 0.077 | 0.271 | 0.368 | 0.277 |
|  | C2 | 1 | 3 | 19 | 43 | 24 |
|  |  | 0.011 | 0.033 | 0.211 | 0.478 | 0.267 |
| I provide time for students to reflect about the processes they use to solve problems: | PC | 1 | 5 | 12 | 28 | 11 |
|  |  | 0.018 | 0.088 | 0.211 | 0.491 | 0.193 |
|  | C1 | 8 | 12 | 59 | 43 | 33 |
|  |  | 0.052 | 0.077 | 0.381 | 0.277 | 0.213 |
|  | C2 | 2 | 16 | 29 | 30 | 14 |
|  |  | 0.022 | 0.176 | 0.319 | 0.33 | 0.154 |
| I give students frequent assignments worth a small portion of their grade: | PC | 6 | 2 | 8 | 15 | 25 |
|  |  | 0.107 | 0.036 | 0.143 | 0.268 | 0.446 |
|  | C1 | 15 | 14 | 36 | 44 | 46 |
|  |  | 0.097 | 0.09 | 0.232 | 0.284 | 0.297 |
|  | C2 | 11 | 8 | 12 | 24 | 36 |
|  |  | 0.121 | 0.088 | 0.132 | 0.264 | 0.396 |
| I expect students to make connections between related ideas or concepts when completing assignments: | PC | 1 | 3 | 13 | 22 | 18 |
|  |  | 0.018 | 0.053 | 0.228 | 0.386 | 0.316 |
|  | C1 | 1 | 13 | 32 | 68 | 41 |
|  |  | 0.006 | 0.084 | 0.206 | 0.439 | 0.265 |
|  | C2 | 0 | 2 | 15 | 46 | 28 |
|  |  | 0 | 0.022 | 0.165 | 0.505 | 0.308 |
| I provide feedback on student assignments without assigning a formal grade: | PC | 21 | 18 | 12 | 5 | 1 |
|  |  | 0.368 | 0.316 | 0.211 | 0.088 | 0.018 |
|  | C1 | 53 | 45 | 25 | 19 | 11 |
|  |  | 0.346 | 0.294 | 0.163 | 0.124 | 0.072 |
|  | C2 | 30 | 24 | 14 | 13 | 10 |
|  |  | 0.33 | 0.264 | 0.154 | 0.143 | 0.11 |

|  |  | Not at all | Minimally | Somewhat | Mostly | Very |
|---|---|---|---|---|---|---|
| Test questions focus on important facts and definitions from the course: | PC | 3 | 7 | 9 | 18 | 20 |
|  |  | 0.053 | 0.123 | 0.158 | 0.316 | 0.351 |
|  | C1 | 3 | 16 | 24 | 57 | 49 |
|  |  | 0.02 | 0.107 | 0.161 | 0.383 | 0.329 |
|  | C2 | 2 | 6 | 25 | 32 | 25 |
|  |  | 0.022 | 0.067 | 0.278 | 0.356 | 0.278 |
| Test questions require students to apply course concepts to unfamiliar situations: | PC | 7 | 18 | 17 | 8 | 6 |
|  |  | 0.125 | 0.321 | 0.304 | 0.143 | 0.107 |
|  | C1 | 14 | 28 | 44 | 34 | 30 |
|  |  | 0.093 | 0.187 | 0.293 | 0.227 | 0.2 |
|  | C2 | 5 | 18 | 26 | 22 | 19 |
|  |  | 0.056 | 0.2 | 0.289 | 0.244 | 0.211 |
| Test questions contain well-defined problems with one correct solution: | PC | 4 | 3 | 10 | 14 | 27 |
|  |  | 0.069 | 0.052 | 0.172 | 0.241 | 0.466 |
|  | C1 | 2 | 3 | 25 | 56 | 64 |
|  |  | 0.013 | 0.02 | 0.167 | 0.373 | 0.427 |
|  | C2 | 1 | 2 | 24 | 32 | 30 |
|  |  | 0.011 | 0.022 | 0.27 | 0.36 | 0.337 |
| I use a grading curve as needed to adjust student scores: | PC | 23 | 11 | 10 | 9 | 4 |
|  |  | 0.404 | 0.193 | 0.175 | 0.158 | 0.07 |
|  | C1 | 54 | 22 | 30 | 19 | 25 |
|  |  | 0.36 | 0.147 | 0.2 | 0.127 | 0.167 |
|  | C2 | 25 | 14 | 21 | 14 | 16 |
|  |  | 0.278 | 0.156 | 0.233 | 0.156 | 0.178 |
| A wide range of students respond to my questions in class: | PC | 2 | 16 | 23 | 11 | 6 |
|  |  | 0.034 | 0.276 | 0.397 | 0.19 | 0.103 |
|  | C1 | 7 | 29 | 61 | 42 | 17 |
|  |  | 0.045 | 0.186 | 0.391 | 0.269 | 0.109 |
|  | C2 | 1 | 13 | 36 | 29 | 12 |
|  |  | 0.011 | 0.143 | 0.396 | 0.319 | 0.132 |

|  |  | Not at all | Minimally | Somewhat | Mostly | Very |
|---|---|---|---|---|---|---|
| I know most of my students by name: | PC | 3 | 5 | 6 | 11 | 32 |
|  |  | 0.053 | 0.088 | 0.105 | 0.193 | 0.561 |
|  | C1 | 4 | 7 | 14 | 15 | 116 |
|  |  | 0.026 | 0.045 | 0.09 | 0.096 | 0.744 |
|  | C2 | 0 | 3 | 8 | 14 | 66 |
|  |  | 0 | 0.033 | 0.088 | 0.154 | 0.725 |
| When calling on students in class, I use randomized response strategies (e.g., picking names from a hat): | PC | 34 | 12 | 7 | 2 | 3 |
|  |  | 0.586 | 0.207 | 0.121 | 0.034 | 0.052 |
|  | C1 | 92 | 27 | 19 | 11 | 7 |
|  |  | 0.59 | 0.173 | 0.122 | 0.071 | 0.045 |
|  | C2 | 53 | 15 | 11 | 8 | 4 |
|  |  | 0.582 | 0.165 | 0.121 | 0.088 | 0.044 |
| I structure class to encourage peer-to-peer support among students | PC | 4 | 9 | 12 | 18 | 13 |
|  |  | 0.071 | 0.161 | 0.214 | 0.321 | 0.232 |
|  | C1 | 17 | 27 | 36 | 30 | 46 |
|  |  | 0.109 | 0.173 | 0.231 | 0.192 | 0.295 |
|  | C2 | 9 | 15 | 25 | 23 | 19 |
|  |  | 0.099 | 0.165 | 0.275 | 0.253 | 0.209 |
| There is a sense of community among the students in my class: | PC | 2 | 10 | 18 | 19 | 7 |
|  |  | 0.036 | 0.179 | 0.321 | 0.339 | 0.125 |
|  | C1 | 6 | 11 | 45 | 61 | 34 |
|  |  | 0.038 | 0.07 | 0.287 | 0.389 | 0.217 |
|  | C2 | 0 | 5 | 36 | 27 | 23 |
|  |  | 0 | 0.055 | 0.396 | 0.297 | 0.253 |
| I require students to work in predetermined or randomized groups: | PC | 22 | 6 | 9 | 7 | 13 |
|  |  | 0.386 | 0.105 | 0.158 | 0.123 | 0.228 |
|  | C1 | 38 | 30 | 26 | 27 | 34 |
|  |  | 0.245 | 0.194 | 0.168 | 0.174 | 0.219 |
|  | C2 | 32 | 15 | 12 | 11 | 21 |
|  |  | 0.352 | 0.165 | 0.132 | 0.121 | 0.231 |

|  |  | Not at all | Minimally | Somewhat | Mostly | Very |
|---|---|---|---|---|---|---|
| I use strategies that have been shown to support students from underrepresented groups: | PC | 11<br>0.193 | 16<br>0.281 | 14<br>0.246 | 9<br>0.158 | 7<br>0.123 |
|  | C1 | 36<br>0.242 | 35<br>0.235 | 44<br>0.295 | 26<br>0.174 | 8<br>0.054 |
|  | C2 | 20<br>0.235 | 25<br>0.294 | 22<br>0.259 | 10<br>0.118 | 8<br>0.094 |
| I consider students' thinking/understanding when planning lessons: | PC | 2<br>0.035 | 2<br>0.035 | 5<br>0.088 | 25<br>0.439 | 23<br>0.404 |
|  | C1 | 3<br>0.019 | 12<br>0.077 | 21<br>0.135 | 51<br>0.329 | 68<br>0.439 |
|  | C2 | 3<br>0.033 | 0<br>0 | 12<br>0.132 | 30<br>0.33 | 46<br>0.505 |
| I use a variety of approaches (e.g., questioning, discussion, formal/informal assessments) to gauge where my students are in their understanding of concepts: | PC | 2<br>0.035 | 1<br>0.018 | 12<br>0.211 | 21<br>0.368 | 21<br>0.368 |
|  | C1 | 5<br>0.032 | 17<br>0.11 | 35<br>0.226 | 50<br>0.323 | 48<br>0.31 |
|  | C2 | 1<br>0.011 | 8<br>0.088 | 33<br>0.363 | 25<br>0.275 | 24<br>0.264 |
| I understand students' previous conceptions, skills, knowledge, and interests related to a particular topic: | PC | 1<br>0.017 | 5<br>0.086 | 19<br>0.328 | 23<br>0.397 | 10<br>0.172 |
|  | C1 | 6<br>0.038 | 16<br>0.103 | 50<br>0.321 | 61<br>0.391 | 23<br>0.147 |
|  | C2 | 2<br>0.022 | 8<br>0.088 | 27<br>0.297 | 39<br>0.429 | 15<br>0.165 |
| I explain concepts in this class in a variety of ways: | PC | 1<br>0.017 | 0<br>0 | 11<br>0.19 | 28<br>0.483 | 18<br>0.31 |
|  | C1 | 1<br>0.006 | 4<br>0.026 | 34<br>0.218 | 69<br>0.442 | 48<br>0.308 |
|  | C2 | 0<br>0 | 2<br>0.022 | 15<br>0.167 | 39<br>0.433 | 34<br>0.378 |

|  |  | Not at all | Minimally | Somewhat | Mostly | Very |
|---|---|---|---|---|---|---|
| I adjust my teaching based upon what students currently do or do not understand: | PC | 2<br>0.035 | 3<br>0.053 | 13<br>0.228 | 17<br>0.298 | 22<br>0.386 |
|  | C1 | 0<br>0 | 6<br>0.038 | 29<br>0.186 | 54<br>0.346 | 67<br>0.429 |
|  | C2 | 1<br>0.011 | 3<br>0.033 | 14<br>0.154 | 39<br>0.429 | 34<br>0.374 |
| I give feedback on homework, exams, quizzes, etc.: | PC | 3<br>0.052 | 3<br>0.052 | 7<br>0.121 | 22<br>0.379 | 23<br>0.397 |
|  | C1 | 4<br>0.026 | 10<br>0.065 | 21<br>0.135 | 40<br>0.258 | 80<br>0.516 |
|  | C2 | 4<br>0.044 | 2<br>0.022 | 9<br>0.099 | 26<br>0.286 | 50<br>0.549 |
| I structure class so that students share their ideas (or their group's ideas) during whole class discussions: | PC | 5<br>0.088 | 5<br>0.088 | 18<br>0.316 | 17<br>0.298 | 12<br>0.211 |
|  | C1 | 15<br>0.096 | 33<br>0.212 | 40<br>0.256 | 47<br>0.301 | 21<br>0.135 |
|  | C2 | 6<br>0.067 | 21<br>0.233 | 21<br>0.233 | 32<br>0.356 | 10<br>0.111 |
| I use strategies to encourage participation from a wide range of students: | PC | 2<br>0.035 | 2<br>0.035 | 20<br>0.351 | 24<br>0.421 | 9<br>0.158 |
|  | C1 | 4<br>0.026 | 22<br>0.145 | 49<br>0.322 | 47<br>0.309 | 30<br>0.197 |
|  | C2 | 3<br>0.033 | 11<br>0.121 | 26<br>0.286 | 33<br>0.363 | 18<br>0.198 |
| A wide range of students participate in class: | PC | 2<br>0.035 | 6<br>0.105 | 25<br>0.439 | 13<br>0.228 | 11<br>0.193 |
|  | C1 | 6<br>0.038 | 22<br>0.14 | 53<br>0.338 | 43<br>0.274 | 33<br>0.21 |
|  | C2 | 0<br>0 | 7<br>0.077 | 22<br>0.242 | 42<br>0.462 | 20<br>0.22 |

## Comparisons

Generally speaking, do other people in your role for [course] use a teaching style similar to yours?

|    | Yes   | No    | Too varied to choose | I don't know |
|----|-------|-------|----------------------|--------------|
| PC | 28    | 4     | 6                    | 20           |
|    | 0.483 | 0.069 | 0.103                | 0.345        |
| C1 | 74    | 4     | 27                   | 55           |
|    | 0.462 | 0.025 | 0.169                | 0.344        |
| C2 | 43    | 4     | 10                   | 36           |
|    | 0.462 | 0.043 | 0.108                | 0.387        |

How do you feel about the instructional approach(es) being used to teach [course] at [your institution]?

|    | Very unhappy | Somewhat unhappy | Neutral | Somewhat happy | Happy |
|----|--------------|------------------|---------|----------------|-------|
| PC | 3            | 1                | 12      | 26             | 17    |
|    | 0.051        | 0.017            | 0.203   | 0.441          | 0.288 |
| C1 | 6            | 10               | 29      | 73             | 42    |
|    | 0.038        | 0.062            | 0.181   | 0.456          | 0.262 |
| C2 | 2            | 6                | 21      | 36             | 27    |
|    | 0.022        | 0.065            | 0.228   | 0.391          | 0.293 |

# PIPS: Instructor Survey

The PIPS data contain 457 responses from instructors from twelve universities across the United States. Of those instructors, 186 taught Precalculus, 164 taught Calculus 1, and 107 taught Calculus 2. We first present the self-identified demographic information of these instructors.

As noted in the introduction to this report, instructors responded for each course they were teaching in a particular term. In the identity and context portion of this report, instructors are counted only once per term; in the course-specific information their responses are reported separately for each *course* they were teaching at the time they completed the survey.

## Instructor Identity & Personal Context

These questions asked about instructors' identity, including demographic information and academic context. Questions about Gender, Race and Ethnicity, and Sexual Orientation were presented as "select all that apply" items with the option to self-identify descriptors not on our lists. We intentionally provided more nuanced options than are generally use in demographic reports (e.g., the US census) as we value the complex nature of identity. Responses to these items were recoded for interpretability in this report; details are provided alongside each item.

### Gender

Do you consider yourself to be (Select all that apply):

| Cisgender Man | Cisgender Woman | Gender Non-binary | Transgender | Transgender Man | Transgender Women | Not reported |
|---|---|---|---|---|---|---|
| 228 | 166 | 5 | 0 | 0 | 0 | 57 |
| 0.5 | 0.364 | 0.011 | 0.000 | 0.000 | 0.000 | 0.125 |

*Survey response options: Man, Transgender, Woman, Not listed (please specify), Prefer not to disclose.*

*Selections of only Man, only Woman, or only Transgender are coded as Cisgender Man, Cisgender Woman, and Transgender respectively; selections of Transgender and Man or Transgender and Woman are coded as Transgender Man and Transgender Woman, respectively; selections of Man and Woman (with or without additional selections) or self-described non-binary gender identities (e.g., genderfluid, agender) are coded as Gender Non-Binary; selections of Prefer not to disclose and no selection were coded as Not reported.*

### Race and Ethnicity

Do you consider yourself to be (Select all that apply):

| | |
|---|---|
| Asian | 37 (0.081) |
| Alaska Native or Native American | 4 (0.009) |
| Black or African American | 43 (0.094) |
| Hispanic or Latinx | 17 (0.037) |
| Middle Eastern or North African | 0 (0.000) |
| Multiple Race | 11 (0.024) |
| Native Hawaiian or Pacific Islander | 0 (0.000) |
| White | 265 (0.58) |
| Not reported | 80 (0.175) |

*Survey response options: Alaskan Native or Native American, Black or African American, Central Asian, East Asian, Hispanic or Latinx, Middle Eastern or North African, Native Hawaiian or Pacific Islander, Southeast Asian, South Asian, White, Not listed (please specify), Prefer not to disclose.*

*Selections of more than one of Central Asian, East Asian, Southeast Asian, South Asian without any other choices are coded as Asian; selections of multiple responses (aside from the previous) are coded as Multiple Race/Ethnicity; selections of Prefer not to disclose and no selection are coded as Not reported.*

### Sexual Orientation

Do you consider yourself to be (Select all that apply):

| Straight | Gay | Lesbian | Bisexual | Asexual | Queer | Queer+ | Straight+ | Not reported |
|---|---|---|---|---|---|---|---|---|
| 321 | 1 | 0 | 2 | 10 | 1 | 1 | 0 | 121 |
| 0.702 | 0.002 | 0.000 | 0.004 | 0.022 | 0.002 | 0.002 | 0.000 | 0.265 |

*Survey response options: Asexual, Bisexual, Gay, Lesbian, Queer, Straight (Heterosexual), Not listed (please specify), Prefer not to disclose.*

*Categories are based on research from Voigt (2020) where selections of exactly one listed option are not recoded; selections of Straight (Heterosexual) and one or more additional options are coded as Straight+; selections of more than one option but not Straight (Heterosexual) are coded as Queer+; selections of Prefer not to disclose and no selection are coded as Not reported.*

### Additional Context / Identity Markers

Do you consider yourself to be (Select all that apply):

| | |
|---|---|
| International Student | 113 (0.247) |
| First-generation college student | 105 (0.23) |
| First-generation higher education | 173 (0.379) |
| Person with a disability | 9 (0.02) |
| English language learner | 77 (0.168) |
| Parent or guardian | 160 (0.35) |

## Course-Specific Information

This section describes the participants perceptions of their role in Precalculus, Calculus I, and/or Calculus II (it may be the case that one participant taught more than one of the above courses during the semester). Each result will be shown in the order Precalculus, Calculus I, Calculus II.

### Course Decisions

How are most decisions about course content (e.g., syllabi, exams, homework, pacing, grading) made for [course]? Clarify if you wish.

|    | I make most decisions | I am part of a team that makes most decisions | Someone else makes most decisions |
|----|-----------------------|-----------------------------------------------|-----------------------------------|
| PC | 68                    | 83                                            | 33                                |
|    | 0.37                  | 0.451                                         | 0.179                             |
| C1 | 71                    | 53                                            | 34                                |
|    | 0.449                 | 0.335                                         | 0.215                             |
| C2 | 50                    | 35                                            | 14                                |
|    | 0.505                 | 0.354                                         | 0.141                             |

How are most decisions about instructional approach (e.g., use of clickers, group work, active learning) made for [course]? Clarify if you wish.

|    | I make most decisions | I am part of a team that makes most decisions | Someone else makes most decisions |
|----|-----------------------|-----------------------------------------------|-----------------------------------|
| PC | 123                   | 51                                            | 11                                |
|    | 0.665                 | 0.276                                         | 0.059                             |
| C1 | 116                   | 30                                            | 14                                |
|    | 0.725                 | 0.188                                         | 0.088                             |
| C2 | 89                    | 9                                             | 1                                 |
|    | 0.899                 | 0.091                                         | 0.01                              |

## Percentages of Class Time

What percent of regular class time in [course], over the whole term, did your students spend… [must total 100]

|  | n | mean | sd | median |
|---|---|---|---|---|
| **All Instructors** | | | | |
| *Working on tasks individually* | 449 | 12.65 | 14.41 | 10 |
| *Working on tasks in small groups* | 449 | 15.65 | 16.50 | 11 |
| *Engaging in whole class discussion* | 449 | 13.00 | 13.30 | 10 |
| *Listening to the instructor lecture or solve problems* | 449 | 58.69 | 22.57 | 60 |
| **Precalculus Instructors** | | | | |
| *Working on tasks individually* | 185 | 16.78 | 17.94 | 11 |
| *Working on tasks in small groups* | 185 | 15.85 | 17.55 | 10 |
| *Engaging in whole class discussion* | 185 | 12.91 | 12.76 | 10 |
| *Listening to the instructor lecture or solve problems* | 185 | 54.46 | 22.15 | 54 |
| **Calculus 1 Instructors** | | | | |
| *Working on tasks individually* | 161 | 9.83 | 10.68 | 8 |
| *Working on tasks in small groups* | 161 | 17.03 | 16.61 | 15 |
| *Engaging in whole class discussion* | 161 | 11.87 | 12.13 | 10 |
| *Listening to the instructor lecture or solve problems* | 161 | 61.27 | 21.77 | 64 |
| **Calculus 2 Instructors** | | | | |
| *Working on tasks individually* | 103 | 9.63 | 10.01 | 8 |
| *Working on tasks in small groups* | 103 | 13.15 | 14.09 | 10 |
| *Engaging in whole class discussion* | 103 | 14.95 | 15.69 | 10 |
| *Listening to the instructor lecture or solve problems* | 103 | 62.27 | 23.51 | 66 |

## Classroom Experience

Please indicate the degree to which the following statements are descriptive of your teaching in [course]:

|  |  | Not at all descriptive | Minimally descriptive | Somewhat descriptive | Mostly descriptive | Very descriptive |
|---|---|---|---|---|---|---|
| I guide students through major topics as they listen: | PC | 2 | 11 | 45 | 64 | 61 |
|  |  | 0.011 | 0.06 | 0.246 | 0.35 | 0.333 |
|  | C1 | 1 | 12 | 32 | 60 | 53 |
|  |  | 0.006 | 0.076 | 0.203 | 0.38 | 0.335 |
|  | C2 | 2 | 4 | 23 | 43 | 29 |
|  |  | 0.02 | 0.04 | 0.228 | 0.426 | 0.287 |

|  |  | Not at all descriptive | Minimally descriptive | Somewhat descriptive | Mostly descriptive | Very descriptive |
|---|---|---|---|---|---|---|
| I provide activities that connect course content to my students' lives and future work: | PC | 17 | 53 | 58 | 28 | 27 |
|  |  | 0.093 | 0.29 | 0.317 | 0.153 | 0.148 |
|  | C1 | 26 | 41 | 51 | 23 | 18 |
|  |  | 0.164 | 0.258 | 0.321 | 0.145 | 0.113 |
|  | C2 | 19 | 26 | 34 | 16 | 5 |
|  |  | 0.19 | 0.26 | 0.34 | 0.16 | 0.05 |
| My syllabus contains the specific topics that will be covered in every class session: | PC | 16 | 11 | 22 | 40 | 94 |
|  |  | 0.087 | 0.06 | 0.12 | 0.219 | 0.514 |
|  | C1 | 27 | 10 | 31 | 33 | 59 |
|  |  | 0.169 | 0.062 | 0.194 | 0.206 | 0.369 |
|  | C2 | 20 | 9 | 11 | 20 | 42 |
|  |  | 0.196 | 0.088 | 0.108 | 0.196 | 0.412 |
| I provide students with immediate feedback on their word during class (e.g., student response systems; short quizzes): | PC | 20 | 26 | 33 | 54 | 50 |
|  |  | 0.109 | 0.142 | 0.18 | 0.295 | 0.273 |
|  | C1 | 40 | 40 | 26 | 26 | 28 |
|  |  | 0.25 | 0.25 | 0.162 | 0.162 | 0.175 |
|  | C2 | 22 | 23 | 21 | 18 | 16 |
|  |  | 0.22 | 0.23 | 0.21 | 0.18 | 0.16 |
| I structure my course with the assumption that most of the students have little useful knowledge of the topics: | PC | 13 | 31 | 42 | 58 | 40 |
|  |  | 0.071 | 0.168 | 0.228 | 0.315 | 0.217 |
|  | C1 | 9 | 24 | 45 | 41 | 39 |
|  |  | 0.057 | 0.152 | 0.285 | 0.259 | 0.247 |
|  | C2 | 6 | 18 | 28 | 35 | 15 |
|  |  | 0.059 | 0.176 | 0.275 | 0.343 | 0.147 |
| I use student assessment results to guide the direction of my instruction during the semester: | PC | 15 | 24 | 50 | 49 | 43 |
|  |  | 0.083 | 0.133 | 0.276 | 0.271 | 0.238 |
|  | C1 | 12 | 22 | 56 | 40 | 29 |
|  |  | 0.075 | 0.138 | 0.352 | 0.252 | 0.182 |
|  | C2 | 5 | 19 | 30 | 24 | 22 |
|  |  | 0.05 | 0.19 | 0.3 | 0.24 | 0.22 |

|  |  | Not at all descriptive | Minimally descriptive | Somewhat descriptive | Mostly descriptive | Very descriptive |
|---|---|---|---|---|---|---|
| I ask students to respond to questions during class time: | PC | 6<br>0.033 | 3<br>0.017 | 5<br>0.028 | 63<br>0.348 | 104<br>0.575 |
|  | C1 | 2<br>0.013 | 2<br>0.013 | 23<br>0.144 | 45<br>0.281 | 88<br>0.55 |
|  | C2 | 1<br>0.01 | 0<br>0 | 12<br>0.121 | 27<br>0.273 | 59<br>0.596 |
| I use student questions and comments to determine the focus and direction of classroom lessons: | PC | 10<br>0.054 | 15<br>0.082 | 58<br>0.315 | 60<br>0.326 | 41<br>0.223 |
|  | C1 | 7<br>0.044 | 21<br>0.131 | 51<br>0.319 | 51<br>0.319 | 30<br>0.188 |
|  | C2 | 1<br>0.01 | 14<br>0.137 | 28<br>0.275 | 40<br>0.392 | 19<br>0.186 |
| In my class a variety of means (models, drawings, graphs, symbols, simulations, tables, etc.) are used to represent course topics and/or solve problems: | PC | 6<br>0.013 | 32<br>0.072 | 98<br>0.22 | 138<br>0.309 | 172<br>0.386 |
|  | C1 | 1<br>0.006 | 11<br>0.069 | 37<br>0.231 | 46<br>0.288 | 65<br>0.406 |
|  | C2 | 0<br>0 | 11<br>0.108 | 28<br>0.275 | 29<br>0.284 | 34<br>0.333 |
| I structure class so that students explore or discuss their understanding of concepts before direct instruction: | PC | 42<br>0.231 | 42<br>0.231 | 43<br>0.236 | 33<br>0.181 | 22<br>0.121 |
|  | C1 | 35<br>0.22 | 53<br>0.333 | 32<br>0.201 | 18<br>0.113 | 21<br>0.132 |
|  | C2 | 28<br>0.286 | 31<br>0.316 | 23<br>0.235 | 8<br>0.082 | 8<br>0.082 |
| My class sessions are structured to give students a clear/structured set of notes: | PC | 1<br>0.005 | 10<br>0.054 | 20<br>0.108 | 52<br>0.281 | 102<br>0.551 |
|  | C1 | 3<br>0.019 | 17<br>0.107 | 15<br>0.094 | 50<br>0.314 | 74<br>0.465 |
|  | C2 | 1<br>0.01 | 7<br>0.069 | 17<br>0.167 | 29<br>0.284 | 48<br>0.471 |

|  |  | Not at all descriptive | Minimally descriptive | Somewhat descriptive | Mostly descriptive | Very descriptive |
|---|---|---|---|---|---|---|
| I structure class so that students talk with one another about course topics: | PC | 14 | 38 | 57 | 43 | 29 |
|  |  | 0.077 | 0.21 | 0.315 | 0.238 | 0.16 |
|  | C1 | 18 | 24 | 48 | 38 | 32 |
|  |  | 0.112 | 0.15 | 0.3 | 0.238 | 0.2 |
|  | C2 | 10 | 21 | 27 | 18 | 24 |
|  |  | 0.1 | 0.21 | 0.27 | 0.18 | 0.24 |
| I structure class so that students constructively criticize one another's ideas: | PC | 62 | 43 | 47 | 18 | 11 |
|  |  | 0.343 | 0.238 | 0.26 | 0.099 | 0.061 |
|  | C1 | 55 | 43 | 30 | 20 | 11 |
|  |  | 0.346 | 0.27 | 0.189 | 0.126 | 0.069 |
|  | C2 | 31 | 29 | 21 | 12 | 6 |
|  |  | 0.313 | 0.293 | 0.212 | 0.121 | 0.061 |
| I structure class so that students discuss their mathematical difficulties with other students: | PC | 24 | 37 | 60 | 40 | 23 |
|  |  | 0.13 | 0.201 | 0.326 | 0.217 | 0.125 |
|  | C1 | 27 | 27 | 51 | 36 | 18 |
|  |  | 0.17 | 0.17 | 0.321 | 0.226 | 0.113 |
|  | C2 | 20 | 14 | 33 | 19 | 14 |
|  |  | 0.2 | 0.14 | 0.33 | 0.19 | 0.14 |
| I structure class so that students work on problems individually during class: | PC | 10 | 39 | 53 | 44 | 37 |
|  |  | 0.055 | 0.213 | 0.29 | 0.24 | 0.202 |
|  | C1 | 27 | 44 | 42 | 23 | 24 |
|  |  | 0.169 | 0.275 | 0.262 | 0.144 | 0.15 |
|  | C2 | 19 | 28 | 28 | 13 | 14 |
|  |  | 0.186 | 0.275 | 0.275 | 0.127 | 0.137 |
| I structure class so that students work together in pairs or small groups: | PC | 29 | 35 | 35 | 47 | 37 |
|  |  | 0.158 | 0.191 | 0.191 | 0.257 | 0.202 |
|  | C1 | 29 | 30 | 37 | 30 | 34 |
|  |  | 0.181 | 0.188 | 0.231 | 0.188 | 0.212 |
|  | C2 | 25 | 23 | 23 | 8 | 22 |
|  |  | 0.248 | 0.228 | 0.228 | 0.079 | 0.218 |

|  |  | Not at all descriptive | Minimally descriptive | Somewhat descriptive | Mostly descriptive | Very descriptive |
|---|---|---|---|---|---|---|
| I structure class so that more than one approach to solving a problem is discussed: | PC | 1<br>0.006 | 7<br>0.039 | 37<br>0.206 | 74<br>0.411 | 61<br>0.339 |
|  | C1 | 1<br>0.006 | 13<br>0.081 | 43<br>0.269 | 70<br>0.438 | 33<br>0.206 |
|  | C2 | 4<br>0.04 | 9<br>0.09 | 28<br>0.28 | 31<br>0.31 | 28<br>0.28 |
| I provide time for students to reflect about the processes they use to solve problems: | PC | 7<br>0.038 | 24<br>0.13 | 60<br>0.326 | 66<br>0.359 | 27<br>0.147 |
|  | C1 | 12<br>0.075 | 26<br>0.164 | 54<br>0.34 | 38<br>0.239 | 29<br>0.182 |
|  | C2 | 4<br>0.04 | 23<br>0.228 | 33<br>0.327 | 26<br>0.257 | 15<br>0.149 |
| I give students frequent assignments worth a small portion of their grade: | PC | 3<br>0.016 | 11<br>0.059 | 28<br>0.151 | 65<br>0.351 | 78<br>0.422 |
|  | C1 | 10<br>0.063 | 8<br>0.05 | 26<br>0.164 | 52<br>0.327 | 63<br>0.396 |
|  | C2 | 2<br>0.02 | 10<br>0.099 | 17<br>0.168 | 29<br>0.287 | 43<br>0.426 |
| I expect students to make connections between related ideas or concepts when completing assignments: | PC | 1<br>0.005 | 4<br>0.022 | 33<br>0.179 | 79<br>0.429 | 67<br>0.364 |
|  | C1 | 3<br>0.019 | 12<br>0.075 | 38<br>0.238 | 71<br>0.444 | 36<br>0.225 |
|  | C2 | 1<br>0.01 | 7<br>0.07 | 23<br>0.23 | 38<br>0.38 | 31<br>0.31 |
| I provide feedback on student assignments without assigning a formal grade: | PC | 65<br>0.357 | 42<br>0.231 | 45<br>0.247 | 18<br>0.099 | 12<br>0.066 |
|  | C1 | 61<br>0.384 | 39<br>0.245 | 38<br>0.239 | 6<br>0.038 | 15<br>0.094 |
|  | C2 | 41<br>0.418 | 22<br>0.224 | 20<br>0.204 | 8<br>0.082 | 7<br>0.071 |

|  |  | Not at all descriptive | Minimally descriptive | Somewhat descriptive | Mostly descriptive | Very descriptive |
|---|---|---|---|---|---|---|
| Test questions focus on important facts and definitions from the course: | PC | 5 | 19 | 33 | 44 | 82 |
|  |  | 0.027 | 0.104 | 0.18 | 0.24 | 0.448 |
|  | C1 | 2 | 12 | 38 | 44 | 63 |
|  |  | 0.013 | 0.075 | 0.239 | 0.277 | 0.396 |
|  | C2 | 0 | 6 | 24 | 31 | 37 |
|  |  | 0 | 0.061 | 0.245 | 0.316 | 0.378 |
| Test questions require students to apply course concepts to unfamiliar situations: | PC | 22 | 50 | 70 | 29 | 14 |
|  |  | 0.119 | 0.27 | 0.378 | 0.157 | 0.076 |
|  | C1 | 21 | 37 | 56 | 22 | 21 |
|  |  | 0.134 | 0.236 | 0.357 | 0.14 | 0.134 |
|  | C2 | 15 | 21 | 39 | 16 | 7 |
|  |  | 0.153 | 0.214 | 0.398 | 0.163 | 0.071 |
| Test questions contain well-defined problems with one correct solution: | PC | 0 | 3 | 17 | 64 | 99 |
|  |  | 0 | 0.016 | 0.093 | 0.35 | 0.541 |
|  | C1 | 1 | 4 | 17 | 60 | 76 |
|  |  | 0.006 | 0.025 | 0.108 | 0.38 | 0.481 |
|  | C2 | 2 | 3 | 19 | 41 | 38 |
|  |  | 0.019 | 0.029 | 0.184 | 0.398 | 0.369 |
| I use a grading curve as needed to adjust student scores: | PC | 69 | 53 | 33 | 14 | 15 |
|  |  | 0.375 | 0.288 | 0.179 | 0.076 | 0.082 |
|  | C1 | 45 | 40 | 27 | 19 | 28 |
|  |  | 0.283 | 0.252 | 0.17 | 0.119 | 0.176 |
|  | C2 | 22 | 22 | 23 | 12 | 20 |
|  |  | 0.222 | 0.222 | 0.232 | 0.121 | 0.202 |
| A wide range of students respond to my questions in class: | PC | 12 | 22 | 61 | 62 | 27 |
|  |  | 0.065 | 0.12 | 0.332 | 0.337 | 0.147 |
|  | C1 | 5 | 26 | 68 | 39 | 22 |
|  |  | 0.031 | 0.162 | 0.425 | 0.244 | 0.138 |
|  | C2 | 3 | 13 | 46 | 23 | 17 |
|  |  | 0.029 | 0.127 | 0.451 | 0.225 | 0.167 |

|  |  | Not at all descriptive | Minimally descriptive | Somewhat descriptive | Mostly descriptive | Very descriptive |
|---|---|---|---|---|---|---|
| I know most of my students by name: | PC | 9 | 23 | 27 | 29 | 95 |
|  |  | 0.049 | 0.126 | 0.148 | 0.158 | 0.519 |
|  | C1 | 16 | 29 | 20 | 20 | 75 |
|  |  | 0.1 | 0.181 | 0.125 | 0.125 | 0.469 |
|  | C2 | 5 | 14 | 12 | 14 | 56 |
|  |  | 0.05 | 0.139 | 0.119 | 0.139 | 0.554 |
| When calling on students in class, I use randomized response strategies (e.g., picking names from a hat): | PC | 105 | 32 | 26 | 8 | 11 |
|  |  | 0.577 | 0.176 | 0.143 | 0.044 | 0.06 |
|  | C1 | 94 | 27 | 19 | 11 | 9 |
|  |  | 0.588 | 0.169 | 0.119 | 0.069 | 0.056 |
|  | C2 | 58 | 15 | 12 | 8 | 5 |
|  |  | 0.592 | 0.153 | 0.122 | 0.082 | 0.051 |
| I structure class to encourage peer-to-peer support among students (e.g., ask peer before you ask me, having group roles, developing a group solution to share, etc.): | PC | 20 | 50 | 37 | 44 | 32 |
|  |  | 0.109 | 0.273 | 0.202 | 0.24 | 0.175 |
|  | C1 | 23 | 35 | 42 | 31 | 29 |
|  |  | 0.144 | 0.219 | 0.262 | 0.194 | 0.181 |
|  | C2 | 21 | 22 | 23 | 18 | 16 |
|  |  | 0.21 | 0.22 | 0.23 | 0.18 | 0.16 |
| There is a sense of community among the students in my class: | PC | 4 | 16 | 64 | 62 | 37 |
|  |  | 0.022 | 0.087 | 0.35 | 0.339 | 0.202 |
|  | C1 | 4 | 15 | 39 | 66 | 35 |
|  |  | 0.025 | 0.094 | 0.245 | 0.415 | 0.22 |
|  | C2 | 1 | 12 | 34 | 30 | 25 |
|  |  | 0.01 | 0.118 | 0.333 | 0.294 | 0.245 |
| I require students to work in predetermined or randomized groups: | PC | 75 | 27 | 32 | 23 | 23 |
|  |  | 0.417 | 0.15 | 0.178 | 0.128 | 0.128 |
|  | C1 | 81 | 16 | 20 | 13 | 30 |
|  |  | 0.506 | 0.1 | 0.125 | 0.081 | 0.188 |
|  | C2 | 56 | 18 | 13 | 7 | 6 |
|  |  | 0.56 | 0.18 | 0.13 | 0.07 | 0.06 |

|  |  | Not at all descriptive | Minimally descriptive | Somewhat descriptive | Mostly descriptive | Very descriptive |
|---|---|---|---|---|---|---|
| I use strategies that have been shown to support students from underrepresented groups: | PC | 34<br>0.192 | 35<br>0.198 | 54<br>0.305 | 36<br>0.203 | 18<br>0.102 |
|  | C1 | 33<br>0.216 | 38<br>0.248 | 39<br>0.255 | 20<br>0.131 | 23<br>0.15 |
|  | C2 | 24<br>0.253 | 29<br>0.305 | 21<br>0.221 | 10<br>0.105 | 11<br>0.116 |
| I consider students' thinking/understanding when planning lessons: | PC | 6<br>0.033 | 6<br>0.033 | 20<br>0.109 | 84<br>0.459 | 67<br>0.366 |
|  | C1 | 0<br>0 | 1<br>0.006 | 27<br>0.172 | 64<br>0.408 | 65<br>0.414 |
|  | C2 | 0<br>0 | 0<br>0 | 16<br>0.158 | 36<br>0.356 | 49<br>0.485 |
| I use a variety of approaches (e.g., questioning, discussion, formal/informal assessments) to gauge where my students are in their understanding of concepts: | PC | 2<br>0.011 | 9<br>0.049 | 45<br>0.246 | 71<br>0.388 | 56<br>0.306 |
|  | C1 | 3<br>0.019 | 17<br>0.106 | 52<br>0.325 | 46<br>0.288 | 42<br>0.262 |
|  | C2 | 1<br>0.01 | 12<br>0.119 | 28<br>0.277 | 35<br>0.347 | 25<br>0.248 |
| I understand students' previous conceptions, skills, knowledge, and interests related to a particular topic: | PC | 4<br>0.022 | 18<br>0.098 | 51<br>0.277 | 63<br>0.342 | 48<br>0.261 |
|  | C1 | 2<br>0.013 | 16<br>0.1 | 61<br>0.381 | 54<br>0.338 | 27<br>0.169 |
|  | C2 | 1<br>0.01 | 10<br>0.098 | 33<br>0.324 | 33<br>0.324 | 25<br>0.245 |
| I explain concepts in this class in a variety of ways: | PC | 2<br>0.011 | 3<br>0.016 | 17<br>0.093 | 83<br>0.456 | 77<br>0.423 |
|  | C1 | 1<br>0.006 | 5<br>0.031 | 34<br>0.214 | 65<br>0.409 | 54<br>0.34 |
|  | C2 | 1<br>0.01 | 4<br>0.04 | 18<br>0.178 | 44<br>0.436 | 34<br>0.337 |

|  |  | Not at all descriptive | Minimally descriptive | Somewhat descriptive | Mostly descriptive | Very descriptive |
|---|---|---|---|---|---|---|
| I adjust my teaching based upon what students currently do or do not understand: | PC | 4 | 9 | 41 | 64 | 65 |
|  |  | 0.022 | 0.049 | 0.224 | 0.35 | 0.355 |
|  | C1 | 0 | 7 | 39 | 63 | 51 |
|  |  | 0 | 0.044 | 0.244 | 0.394 | 0.319 |
|  | C2 | 0 | 4 | 28 | 36 | 34 |
|  |  | 0 | 0.039 | 0.275 | 0.353 | 0.333 |
| I give feedback on homework, exams, quizzes, etc.: | PC | 4 | 8 | 26 | 59 | 87 |
|  |  | 0.022 | 0.043 | 0.141 | 0.321 | 0.473 |
|  | C1 | 7 | 10 | 28 | 41 | 73 |
|  |  | 0.044 | 0.063 | 0.176 | 0.258 | 0.459 |
|  | C2 | 4 | 4 | 9 | 30 | 55 |
|  |  | 0.039 | 0.039 | 0.088 | 0.294 | 0.539 |
| I structure class so that students share their ideas (or their group's ideas) during whole class discussions: | PC | 17 | 45 | 50 | 46 | 26 |
|  |  | 0.092 | 0.245 | 0.272 | 0.25 | 0.141 |
|  | C1 | 26 | 32 | 47 | 29 | 26 |
|  |  | 0.162 | 0.2 | 0.294 | 0.181 | 0.162 |
|  | C2 | 10 | 26 | 32 | 20 | 13 |
|  |  | 0.099 | 0.257 | 0.317 | 0.198 | 0.129 |
| I use strategies to encourage participation from a wide range of students: | PC | 7 | 11 | 48 | 69 | 50 |
|  |  | 0.038 | 0.059 | 0.259 | 0.373 | 0.27 |
|  | C1 | 8 | 28 | 38 | 44 | 41 |
|  |  | 0.05 | 0.176 | 0.239 | 0.277 | 0.258 |
|  | C2 | 5 | 21 | 28 | 28 | 17 |
|  |  | 0.051 | 0.212 | 0.283 | 0.283 | 0.172 |
| A wide range of students participate in class: | PC | 7 | 13 | 58 | 72 | 33 |
|  |  | 0.038 | 0.071 | 0.317 | 0.393 | 0.18 |
|  | C1 | 4 | 15 | 60 | 49 | 32 |
|  |  | 0.025 | 0.094 | 0.375 | 0.306 | 0.2 |
|  | C2 | 1 | 15 | 33 | 35 | 18 |
|  |  | 0.01 | 0.147 | 0.324 | 0.343 | 0.176 |

# Comparisons

Generally speaking, do other [course] instructors use a teaching style similar to yours?

|  |  | Yes | No | Too varied to choose | I don't know |
|---|---|---|---|---|---|
| Do other [course] instructors use a teaching style similar to yours? | PC | 72 | 13 | 22 | 77 |
|  |  | 0.391 | 0.071 | 0.12 | 0.418 |
|  | C1 | 57 | 18 | 18 | 64 |
|  |  | 0.363 | 0.115 | 0.115 | 0.408 |
|  | C2 | 30 | 14 | 12 | 45 |
|  |  | 0.297 | 0.139 | 0.119 | 0.446 |

How do you feel about the instructional approach(es) being used to teach [course] at [your institution]?

|  |  | Very unhappy | Somewhat unhappy | Neutral | Somewhat happy | Happy |
|---|---|---|---|---|---|---|
| How do you feel about the instructional approach(es) being used to teach [course] at [your institution]? | PC | 2 | 15 | 40 | 70 | 56 |
|  |  | 0.011 | 0.082 | 0.219 | 0.383 | 0.306 |
|  | C1 | 1 | 26 | 33 | 65 | 34 |
|  |  | 0.006 | 0.164 | 0.208 | 0.409 | 0.214 |
|  | C2 | 0 | 13 | 28 | 34 | 25 |
|  |  | 0 | 0.13 | 0.28 | 0.34 | 0.25 |

# Acknowledgement & References

## Acknowledgements

We are grateful to the students and instructors across all twelve universities who provided us with descriptions of themselves and their mathematics course experiences.

The authors of this manuscript are a subset of the Progress through Calculus research team. This work would not be possible without the contributions of the PI team (David Bressoud, Jessica Ellis Hagman, Sean Larsen, Chris Rasmussen), the Mathematical Association of America, and a long list of senior personnel and graduate students, including: Estrella Johnson, Naneh Apkarian, Jessica Gehrtz, Kristen Vroom, Dana Kirin, Matthew Voigt, Tenchita Alzaga Elizondo, Brittney Ellis, Antonio Martinez, Rachel Tremaine, Ciera Street, Tyler Sullivan, Jason Guglielmo, Victoria Barron, Colin McGrane, and Kate Yang.

This material is based upon work supported by the National Science Foundation under grant DUE IUSE #1430540. Any opinions, findings, and conclusions or recommendations expressed in this material are those of the authors and project team and do not necessarily reflect the views of the National Science Foundation.

## References

Apkarian, N., & Kirin, D. (2017). *Progress through calculus: Census survey technical report*. Mathematical Association of America. http://bit.ly/PtC_Reporting

Apkarian, N., Smith, W. M., Vroom, K., Voigt, M., Gehrtz, J., PtC Project Team, & SEMINAL Project Team. (2019). *X-PIPS-M Survey Suite*. https://www.maa.org/sites/default/files/XPIPSM%20Summary%20Document.pdf

Bressoud, D., Mesa, V., & Rasmussen, C. (Eds.). (2015). *Insights and recommendations from the MAA national study of college calculus*. MAA Press.

Rasmussen, C., Apkarian, N., Hagman, J. E., Johnson, E., Larsen, S., Bressoud, D., & Progress through Calculus Team. (2019). Characteristics of Precalculus through Calculus 2 programs: Insights from a national census survey. *Journal for Research in Mathematics Education*, *50*(1), 98–112.

Voigt, M. (2020). Queer-Spectrum Student Experiences and Resources in Undergraduate Mathematics [Doctoral dissertation, UC San Diego]. Retrieved from https://escholarship.org/uc/item/7g54x6c7

Walter, E. M., Henderson, C. R., Beach, A. L., & Williams, C. T. (2016). Introducing the Postsecondary Instructional Practices Survey (PIPS): A concise, interdisciplinary, and easy-to-score survey. *CBE—Life Sciences Education*, *15*(4), ar53. https://doi.org/10.1187/cbe.15-09-0193